\documentclass[11pt, a4paper]{article}
\usepackage{amsmath, amssymb, amsfonts, amsthm}
\usepackage[pdftex]{color}
\usepackage{subcaption}
\usepackage{graphicx}
\usepackage{float}

\newcommand{\E}{{\bf{E}}}
\newcommand{\PP}{{\bf{P}}}
\newcommand{\Var}{{\bf{Var}}}
\newcommand{\Cov}{{\bf{Cov}}}

\newtheorem{tm}{Theorem}
\newtheorem{lem}{Lemma}
\newtheorem{cor}{Corollary}

\begin{document}

\parindent=0pt

\smallskip
\par\vskip 3.5em
\centerline{\Large \bf  Two models of sparse and  clustered  dynamic  networks}

\vglue2truecm

\vglue1truecm
\centerline {Mindaugas Bloznelis, Dominykas Marma}

\bigskip

\centerline{Institute of Computer Science, Vilnius University}
\centerline{
\, \ Didlaukio 47, LT-08303 Vilnius, Lithuania} 

\vglue2truecm

\abstract{We present two  models of  sparse  dynamic networks
that display transitivity - the tendency for vertices sharing a common neighbour to be neighbours of one another.
Our first  network is a continuous time Markov chain
$G=\{G_t=(V,E_t), t\ge 0\}$ whose states are graphs with the common vertex set $V=\{1,\dots, n\}$. The transitions are defined as follows. Given $t$, the 
 vertex 
pairs $\{i,j\}\subset V$ 
are assigned independent exponential waiting times 
$A_{ij}$. At time $t+\min_{ij} A_{ij}$ the pair $\{i_0,j_0\}$ with 
$A_{i_0j_0}=\min_{ij} A_{ij}$  toggles  its adjacency status.
To mimic clustering patterns of sparse real networks 
we set  intensities $a_{ij}$ of exponential times $A_{ij}$  to be  negatively correlated 
with the degrees of the common neighbours of  vertices $i$ and $j$ in $G_t$.
Another dynamic network is based on a latent Markov chain 
$H=\{H_t=(V\cup W, E_t), t\ge 0\}$ whose states are bipartite graphs with the bipartition $V\cup W$,
where
$W=\{1,\dots,m\}$ is an 
auxiliary 
set of attributes/affiliations. 
Our second network
$G'=\{G'_t =(E'_t,V), t\ge 0\}$ is the affiliation network
defined by  $H$:   vertices $i_1,i_2\in V$  are adjacent in $G'_t$ whenever $i_1$ and $i_2$ 
have  a common neighbour
in $H_t$.  We analyze geometric properties of both dynamic networks at stationarity
and show that networks possess high clustering. They admit tunable degree distribution and clustering coefficients.}
\smallskip
\par\vskip 3.5em

\section{Introduction} Many real networks, especially those depicting human interaction, like  social networks of friendships, collaboration networks,  citation networks and other show  clustering, the propensity of nodes to cluster together by forming
relatively small groups with a high density of ties within a group.  
Clustering is closely related to network transitivity,
the tendency for
two nodes sharing a common neighbour  to be neighbors of one another thus forming a triangle of connections.
Locally, in a vicinity of a node, this tendency can be quantified by the probability 
 that two randomly selected neighbours of the node are adjacent. The network average of this probability, called the (average) local clustering coefficient, is used to quantify the network transitivity.  Another popular measure
 of network transitivity, the  global clustering coefficient,  is  the probability   that two randomly selected neighbours of a randomly selected node are adjacent. 
 In many social networks  both clustering coefficients are on the 
 order of tens of persent while  the edge density, the probability that two randomly selected nodes are adjacent, is of much smaller order.
 Often the edge density  scales as $n^{-1}$, where $n$ is the number of nodes in the network. We call networks with such edge densities sparse.

 Mathematical modelling of sparse networks  displaying clustering/transitivity has attracted considerable attention in the literature, see e.g., \cite{Hofstad_book_2_2024} and references therein. We briefly review several  approaches to modeling of clustered networks. 
 In order to enhance the number of triangles in an evolving  locally tree-like network Holme and Kim \cite{Holme_Kim_2002} suggested inserting additional edges that close desired fraction of open triangles (paths of lenght two).  Newman \cite{MEJNewman_2009} generalised  the configuration random graph model by prescribing network nodes  numbers of triangles  they participate in. In this way a predefined number of triangles can be introduced  into configuration random graph. 
 Bollob\'as et al.
 %
 \cite{Bollobas_Janson_Riordan_2011} built a clustered network by taking a union of randomly located small dense subgraphs of variable sizes. 
Guillaume and Latapy \cite{Guillaume_Latapy_2004} noted  an underlying  bipartite structure present in many social networks, where nodes (actors)  sharing a common hobby or affiliation are more likely to become friends, and where each hobby/affiliation defines a tightly connected cluster of actors related to it. They suggested modelling 
 a clustered network by first linking actors to affiliations
 and  then connecting  actors that share common affiliations, see also \cite{B_G_J_K_R_2015_models}, \cite{Hofstad_Komjathy_Vadon_2021}.  We call such networks  affiliation networks.
 
 The present paper is devoted to the modelling of sparse and clustered dynamic networks 
 using Markov chains. By dynamic network we mean a collection of 
 random graphs  $\{G_t=(V,E_t)$, $t\ge 0\}$ sharing the same vertex 
 set $V=\{1,\dots, n\}$ and having  random edge sets 
 $E_t$, $t\ge 0$.
 We present two stationary random processes $\{G_t, t\ge 0\}$ with tunable 
 degree distribution and tunable non-vanishing clustering coefficients.
 Our study is build upon earlier work on dynamic network 
 Markov chains  \cite{Grindrod_Higham_Parsons_2012}, \cite{UzupyteWit_2020}, 
 \cite{Zhang_Moore_Newman_2017}.
 We mention that  network  Markov chain  of \cite{Zhang_Moore_Newman_2017}
is composed of  ${n}\choose{2}$
independent Markov chains  defining the adjacency status of each 
vertex pair $\{i, j\}\subset V$ individually (we refer to Section 2 for details). The network   admits tunable edge density and degree distribution, but since the edges are inserted/deleted 
 independently of each other it
does not show clustering.
Grindrod et al \cite{Grindrod_Higham_Parsons_2012} 
and Užupytė and Witt \cite{UzupyteWit_2020} introduced  transitivity into the network Markov chain
by relating the birth/death rate of an edge to the number of 
triangles it participates in (cf.  \cite{Holme_Kim_2002}). 
More preciselly, they set the birth (death) 
rate of an edge $i\sim j$ to be an affine function of the number of  
the common neighbours of vertices $i$ and $j$. Here   $i\sim j$
means that $i$ and $j$ are adjacent. 
A drawback of the  models of
 \cite{Grindrod_Higham_Parsons_2012}, \cite{UzupyteWit_2020}
is that for large $n$ they have a little control  over the edge density 
and  clustering strength.
  
  In  the present paper we  suggest a remedy to this drawback.
Inspired by  clustering patterns observed in  real networks, where the 
number of closed triangles incident to a vertex negatively correlates 
with the degree of the vertex (\cite{Foudalis_J_P_Sideri_2011},
\cite{Ravasz_S_M_O_Barabasi_2002},
\cite{Ravasz_Barabasi_2003},  
\cite{Vazquez_Pastor-Satorras_Vespignani_2002})
we set the birth  rate of an edge 
$i\sim j$ to be negatively correlated with the degrees of the common 
neighbours of $i$ and $j$. We show below that such a modification leads to a stationary dynamic 
network model admiting  tunable edge density and clustering 
coefficients.

Another dynamic clustered 
network considered in this paper is a stationary affiliation network  built upon an underlying bipartite graph valued Markov chain
with independent edges.  Now  the clustering property is caused by the bipartite structure  as  noted in \cite{Guillaume_Latapy_2004}.  We analyse the degree sequence and global clustering coefficient at stationarity using the tools developed for random intersection graphs \cite{B_G_J_K_R_2015_models}.
 We note that  earlier work on dynamic affiliation network models 
(\cite{Bloznelis_Goetze_2014}, \cite{Bloznelis_Karonski_2015},
 \cite {Guillaume_Latapy_2006}) addresses the case where the network size $n=n(t)$ increases with time. Clearly, such networks do not admit stationary distributions. 

Finally, we mention  the  recent work by Milewska et al. \cite{Milewska_Hofstad_Zwart}, where a sparse and clustered dynamic network is constructed by taking  unions of small dense subgraphs that are inserted/deleted at random times (cf.  \cite{Bollobas_Janson_Riordan_2011}).

The rest of the paper is organized as follows. In section 2  we formally define 
the network Markov chain and  analyze geometric properties of  the network analytically and by numerical simulations.
In section 3 we define  stationary affiliation network and  show the degree distribution and global clustering coefficient. Proofs of the results of section 3 are given in Appendix. 

\section{Network Markov chain}

Let
${\bf G}=\{G_t=(V,E_t),\, t\ge 0\}$ be a continuous time Markov 
chain, whose states are graphs on the vertex set $V$ and  transitions 
are defined as follows.
Given $G_t$ (the state occupied at time $t$), the update takes place  at time  $t':=t+\min_{ij} A_{ij}$, where $A_{ij}=A_{ij}(G_t)$, 
$\{i,j\}\in V$ are independent exponential waiting times with intensities $a_{ij}=a_{ij}(G_t)$ defined below.  The pair $\{i_0,j_0\}$ with 
$A_{i_0j_0}=\min_{ij} A_{ij}$  changes its adjacency status:  the edge $i_0\sim j_0$ is 
inserted if 
it is not present at time $t$; the edge $i_0\sim j_0$ is 
removed if it is present at time $t$. Thus, at time $t'$ the Markov chain
jumps to the state $G_{t'}=(V, E_{t'})$, where the edge sets $E_t$ and $E_{t'}$ differ in the single edge $i_0\sim j_0$. 

Let us define the intensities $a_{ij}$ for $\{i,j\}\subset V$.
Let   $\alpha, \beta, \lambda,\mu\ge 0$ and let  $\lambda_i, \mu_i$, $1\le i\le n$,
be positive numbers.
  Given graph $G=(V,E)$  we assign clustering weights 
$\nu_{ij}(G,\alpha)$ and $\nu_{ij}(G,\beta)$ to  each vertex pair $\{i,j\}\subset V$, where
\begin{align}
\label{clustering}
\nu_{ij}(G,s)=\sum_{v\in N_{ij}}(d_v(G))^{-s},
\qquad 
s\ge 0.
\end{align}
For $s=0$ we have $\nu_{ij}(G,0)
=|N_{ij}|$.
Here $N_{ij}=N_{ij}(G)$ stands for the set of common neighbours of 
vertices $i$ and $j$ in $G$; $d_v(G)$ denotes the degree of vertex $v$ in 
$G$. 
Furthermore,  each vertex pair $\{i,j\}$ is assigned intensity 
\begin{align}\label{model1}
a_{ij}(G) =
\begin{cases}
\lambda_i\lambda_j+\lambda\nu_{ij}(G, \alpha)
\qquad
\quad
\
\
 {\text{for}}
\quad 
\{i,j\}\not\in E,
\\
\left(\mu_i\mu_j-\mu\nu_{ij}(G,\beta)\right)_+
\qquad
{\text{for}}
\quad
\{i,j\}\in E.
\end{cases}
\end{align}
Here $x_+$ stands for $\max\{x,0\}$. A standard argument shows 
that the chain $\bf G$ has 
 unique stationary distribution. Chain $\bf G$ starting with random graph $G_0$  having such a distribution is called stationary network in what follows.

For $\lambda=\mu=0$  transitions of the chain $\bf G$ are defined 
by the  transitions of
${n \choose 2}$ independent Markov chains describing  adjacency  
dynamic 
of
each vertex pair  $\{i,j\}\subset V$ separately.  (The Markov chain of 
the vertex pair $\{i,j\}$ has two states $i\sim j$ and $i\not\sim j$, 
where 
state $i\sim j$ ($i$ and $j$ are adjacent) has exponential holding time
 with the intensity $\mu_i\mu_j$ and the state $i\not\sim j$ ($i$ and 
 $j$ 
 aren't adjacent) has exponential holding time with the 
 intensity $\lambda_i\lambda_j$.)   The stationary
network of $\bf G$  has independent edges and, hence, it lacks the clustering 
property. Assuming, in addition, that $\mu_i$  is the 
same for each vertex $i\in V$ ($\mu_i\equiv const$) we obtain a dynamic network considered in 
\cite{Zhang_Moore_Newman_2017}. Let us mention that  weights 
$\lambda_i$ strongly correlate with respective vertex degrees $d_i(G_t)$, 
$i\in V$, and  are useful in modeling the degree distribution of $G_t$ for large $t$.
Furthermore, large values of $\lambda_i, \mu_i$  enhance the variability (over time) of links incident to vertex $i\in V$. 

Grindrod et al. \cite{Grindrod_Higham_Parsons_2012} 
introduced the  term  $ \lambda\nu_{ij}(G,0)$ to enhance the triadic closure effect.
We mention that
 \cite{Grindrod_Higham_Parsons_2012} considers the (discrete) jump chain ${\bf G}^*=\{G^*_k=(V,E^*_k), k=0,1,2,\dots\}$ related to    $\bf G$ defined by (\ref{model1}), where    
 $\lambda_i=const_1$,
$\mu_i=const_2$ do not depend on $i$ and where $\mu=0$. More precisely, ${\bf G}^*$  represents the list of  distinct states visited  by the chain
${\bf G}$ arranged  in the chronological order. That is, 
$G^*_0=G_0$, $G^*_1=G_{t_1}$, $G^*_2=G_{t_2}$, $\dots$, where $t_1<t_2<\dots$ are the subsequent jump times  of continuous chain $\bf G$.  
Užupytė and Wit 
\cite{UzupyteWit_2020} complemented the model of  
\cite{Grindrod_Higham_Parsons_2012} by adding the ``triadic 
protection'' term $ \mu\nu_{ij}(G,0)$ aimed at reducing the deletion 
rate of the edges belonging to the closed triangles. They consider
the continuos chain $\bf G$ defined by (\ref{model1}) with $\lambda_i=const_1,
\mu_i=const_2$.

It has already been mentioned that   for large $n$ dynamic networks  of \cite{Grindrod_Higham_Parsons_2012}, \cite{UzupyteWit_2020}
permit little control  over the edge density, which becomes very sensitive to parameters $\mu$ and $\lambda$.
To overcome  such 
disadvantage we suggest choosing clustering weights $\nu_{ij}(G,s)$
that correlate negatively with degrees of the 
common neigbours of $i$ and $j$. An intuition behind this choice is based on the plausible 
assumption that for $i, j$  being  friends of an individual 
with a  large number of acquaintances 
makes less impact on the  mutual relations between $i,j$ than beying friends with 
a person having just a few contacts.
Moreover, \cite{Ravasz_S_M_O_Barabasi_2002},
\cite{Ravasz_Barabasi_2003},  
\cite{Vazquez_Pastor-Satorras_Vespignani_2002}, see also \cite{Foudalis_J_P_Sideri_2011}, note that in some sparse and clustered real networks  the fraction of closed triangles incident to a vertex scales as  a negative power of the  degree of that vertex.
Findings of  \cite{Foudalis_J_P_Sideri_2011},
\cite{Ravasz_S_M_O_Barabasi_2002},
\cite{Ravasz_Barabasi_2003},  
\cite{Vazquez_Pastor-Satorras_Vespignani_2002} motivated our choice of the clustering weights~(\ref{clustering}).

 We are most interested in sparse networks, where the number of vertices $n$ is large. In the simplest case, where 
 $\lambda=\mu=0$ 
 and wehere 
 $\lambda_i=const_1$ and $\mu_i=const_2$ are the same for each $i$ (we write, for short, $\lambda_i\lambda_j=\lambda_0$
 and $\mu_i\mu_j=\mu_0$) 
 each vertex pair toggles its adjacency status independently and  the  expected holding time of an edge (respectively, non-edge)  is $\mu_0^{-1}$ (respectively,  $\lambda_0^{-1}$).
By the law of large numbers the probability that $i$ and $j$ are adjacent in $G_t$  is asymptotically $\mu_0^{-1}/(\mu_0^{-1}+\lambda_0^{-1})=\lambda_0/(\mu_0+\lambda_0)$ as  $t\to\infty$.
 Hence  a snapshot $G_t$ of the stationary network has the distribution of  the binomial random graph with the edge density 
 $\lambda_0/(\mu_0+\lambda_0)$.
 Furthermore, a sparse network is obtained if one chooses $\mu_0=n$ and  $\lambda_0=c$, where $c>0$ denotes a number independent of $n$
(think of a sequence of network Markov chains with vertex number $n\to\infty$).  More generally, for $\lambda=\mu=0$, 
$\mu_i\mu_j\equiv n$ 
 and  $\sum_{i=1}^n\lambda_i\le cn$ uniformly in $n$ one can obtain  a sparse stationary network having independent edges and the degree sequence  strongly correlated with the sequence of weights $\{\lambda_i\}$, \cite{Zhang_Moore_Newman_2017}.

The  simulation study  of subsection 2.1 below shows that  network 
Markov chain (\ref{model1}) with clustering weights
 $\nu_{ij}(G, \alpha)$, 
$\nu_{ij}(G, \beta)$, where $\alpha,\beta>0$,
 can produce  highly clustered sparse stationary dynamic networks with tunable edge density and clustering coefficients.
These empirical findings are supported by a limited analytical study  (given in subsection 2.2 below) showing upper and lower bounds of the order $n^{-1}$ on the average edge density. In addition, we establish a lower bound of the order $n$ on  the average number of triangles
and in a special case of $\alpha=2$ we relate the average edge density to the average local clustering coefficient. 
  
Before proceeding further, we introduce some notation. 
We use terms vertex and node interchangeably.
Given a graph $G=(V,E)$ we denote by $\Delta_v(G)$ the number of triangles incident to a vertex $v\in V$. The total number of triangles is 
denoted $N_{\Delta}(G)=\frac{1}{3}\sum_{v\in V}\Delta_v(G)$
The total number of $2$-paths is 
denoted $N_{\Lambda}(G)=\sum_{v\in V}\binom{d_v(G)}{2}$.
For a vertex $v\in V$ of degree $d_v(G)\ge 2$ we denote 
$C^{\text{\tiny L}}_v(G)=\Delta_v(G){{d_v(G)}\choose{2}}^{-1}$ the 
local clustering coefficient of $v$ ($=$ probability that 
two randomly selected neighbours of $v$ are neighbours to each other). 
In the case where  $d_v(G)\le 1$ we put $C^{\text{\tiny L}}_v(G)=0$.
The average local clustering coefficient and the global clustering coefficient are  denoted 
\begin{displaymath}
{\bar C}^{\text{\tiny L}}(G)
=
\frac{1}{n}\sum_{v\in V}C^{\text{\tiny L}}_v(G)
\qquad
{\text{and}}
\qquad
C^{\text{\tiny GL}}(G)
=
\frac{3N_{\Delta}(G)}
{N_{\Lambda}(G)}.
\end{displaymath}
We put $C^{\text{\tiny GL}}(G)=0$ when 
 $N_{\Delta}(G)=0$. The average degree and the  average edge density are denoted 
 ${\bar d}(G)=n^{-1}\sum_{v\in V}d_v(G)$ 
 and $e(G)={{n}\choose{2}}^{-1}|E|$ respectively. 
Finally, we denote by
${\mathbb I}_{A}$ the indicator function of an event (or set) $A$.

 \subsection{Numerical Simulations}

The aim of the simulation study is twofold: testing the clustering properties of sparse network  (\ref{model1}) equipped with clustering weights $\nu_{ij}(G, \alpha)$, 
$\nu_{ij}(G, \beta)$, where $\alpha,\beta>0$ and comparison of  the clustering properties  for  $\alpha,\beta>0$ and $\alpha=\beta=0$ (the case $\alpha=\beta=0$ corresponds to the setup of
\cite{Grindrod_Higham_Parsons_2012}, 
\cite{UzupyteWit_2020}).

 To address both questions simultaneously 
we consider a simplified model (\ref{model1}), where we assume that  
$\lambda_i\lambda_j\equiv const_1:=\lambda_0$ and $\mu_i\mu_j\equiv const_2:=\mu_0$, see (\ref{modelUW}) below.  Recall that for  $\lambda=\mu=0$ the edges are inserted/deleted independently of each other and   the ratio $\mu_0/\lambda_0$ defines the network edge density $1/(1+\mu_0/\lambda_0)$  at stationarity. Hence, tuning the ratio $\mu_0/\lambda_0$ one can achieve the desired edge density. Here we assume that the ratio $\mu_0/\lambda_0$ is fixed and address the question  about tuning  parameters $\lambda$ and $\mu$ for achieving desired values of clustering coefficients.
 
 In the simulations we put the vertex number  $n=1000$, $\mu_0=n$ and $\lambda_0=1$  (for $\lambda=\mu=0$ such network is sparse at stationarity).
We only consider two instances of values of the pair $(\alpha, \beta)$: the choice 
of parameters
$\alpha=2.75$ and $\beta=2.5$ is referred to as ``general triadic model'' below;
the choice of parameters $\alpha=\beta=0$ is  referred to as  ``simple triadic model''. Given $(\alpha, \beta)$ we generate network Markov chains for different values of $(\mu,\lambda)$ from the range 
that features variability of the 
 local clustering coefficient (our target parameter).
For each choice of  $(\mu,\lambda)$ we sample  network snaphot $G_t$ out of (approximately) stationary distribution  and evaluate the edge density $e(G_t)$ (Figure 1) and local clustering coefficient 
${\bar C}^{\text{\tiny L}}(G_t)$ (Figure 2).
To generate an approximately stationary network we  run  the respective Markov chain starting from 
an empty graph  until $3n^2$ jumps (edge changes) occur. Further simulation steps
 do not 
change values of $e(G_t)$
 and ${\bar C}^{\text{\tiny L}}(G_t)$ beyond the rounding error.  
\begin{figure}[H]
	\centering
	\begin{subfigure}[b]{0.45\textwidth}
		\centering
		\includegraphics[width=\textwidth]{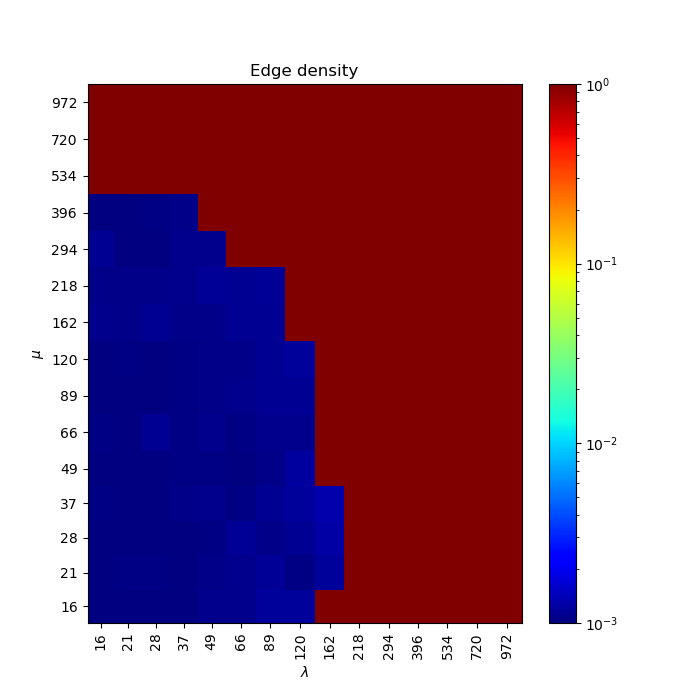} 
		\caption{Simple triadic model}
		\label{fig:Uzupytes_edge_density}
	\end{subfigure}
	\hfill
	\begin{subfigure}[b]{0.45\textwidth}
		\centering
		\includegraphics[width=\textwidth]{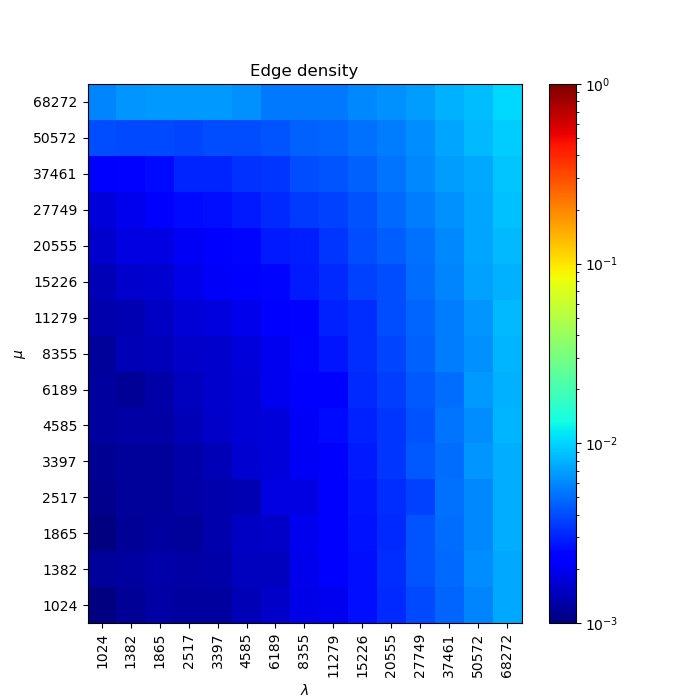}
		\caption{General triadic model}
		\label{fig:alpha_beta_edge_density}
	\end{subfigure}
	\caption{Edge densities in stationary graphs}
	\label{fig:edge_densities}
\end{figure}

\begin{figure}[H]
	\centering
	\begin{subfigure}[b]{0.45\textwidth}
		\centering
		\includegraphics[width=\textwidth]{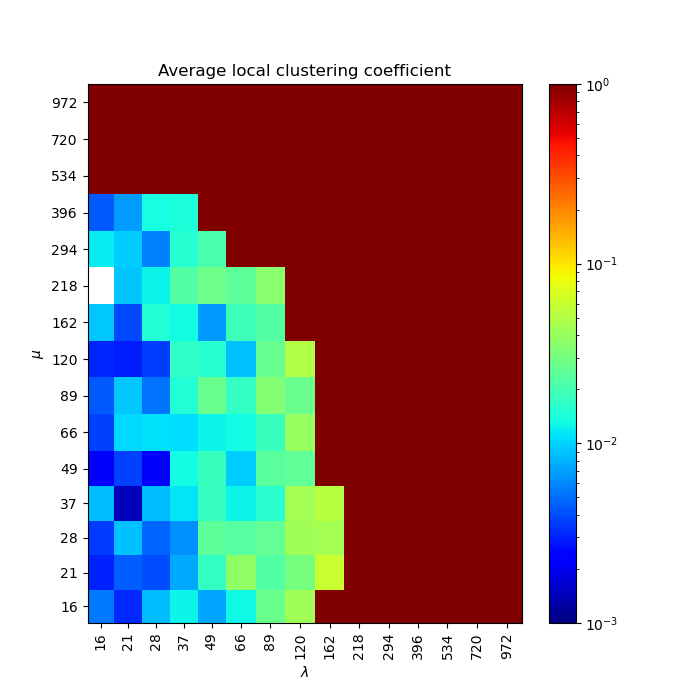}
		\caption{Simple triadic model}
		\label{fig:Uzupytes_local_clustering}
	\end{subfigure}
	\hfill
	\begin{subfigure}[b]{0.45\textwidth}
		\centering
		\includegraphics[width=\textwidth]{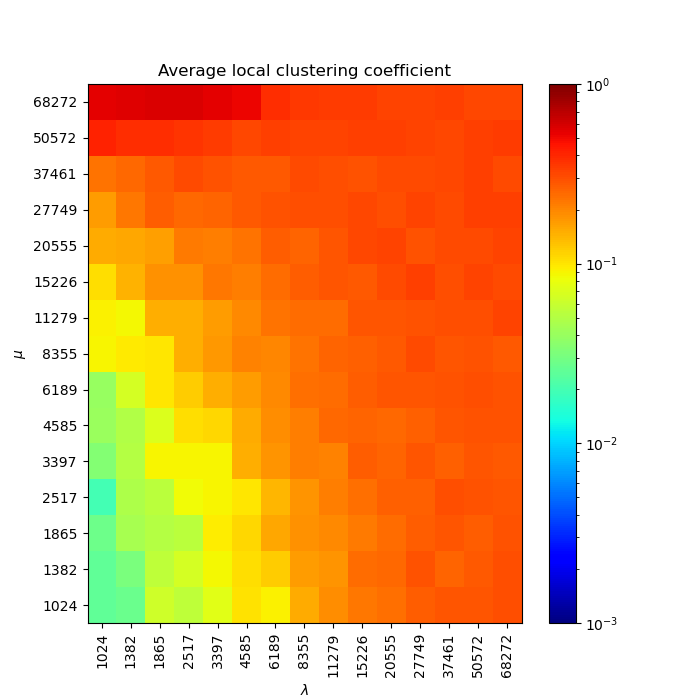}
		\caption{General triadic model}
		\label{fig:alpha_beta_local_clustering}
	\end{subfigure}
	\caption{Average local clustering coefficients in stationary graphs}
	\label{fig:local_clustering}
\end{figure}
\noindent In Figures \ref{fig:edge_densities} and \ref{fig:local_clustering} values of parameters  $\mu$ and $\lambda$ are depicted on the vertical and horizontal axis respectively. Evenly spaced labels on each  axis depict  values of  geometric sequences with the  common ratio $1.35$.
The  colours are put on logarithmic scale and the same scale is applied across different images.


As we can see from Figure \ref{fig:edge_densities}, ``general triadic model'' admits  tunable 
(average) local clustering coefficient while the edge density remains reasonably small (recall that  the ratio $\mu_0/\lambda_0$ remains fixed).
On the other hand, ``simple triadic model'' shows a swift jump from a sparse graph to the complete graph. Hence while trying to achieve the desired values of the clustering coefficient we are losing cotrol over the edge density.

\begin{figure}[H]
	\centering
	\begin{subfigure}[b]{0.45\textwidth}
		\centering
		\includegraphics[width=\textwidth]{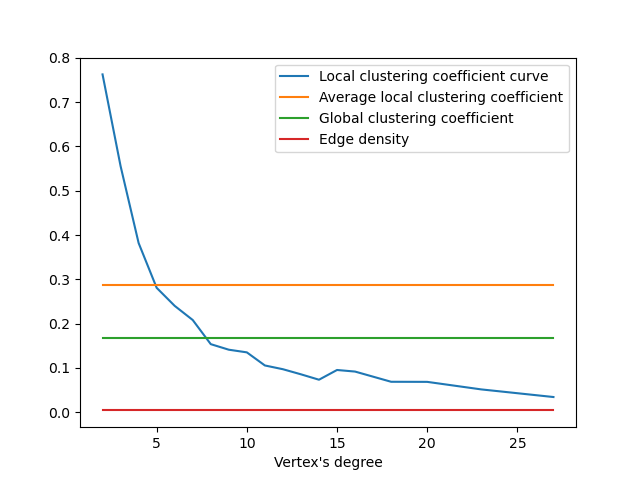}
		\caption{Clustering versus degree}
		\label{fig:Uzupytes_local_clustering}
	\end{subfigure}
	\hfill
	\begin{subfigure}[b]{0.45\textwidth}
		\centering
		\includegraphics[width=\textwidth]{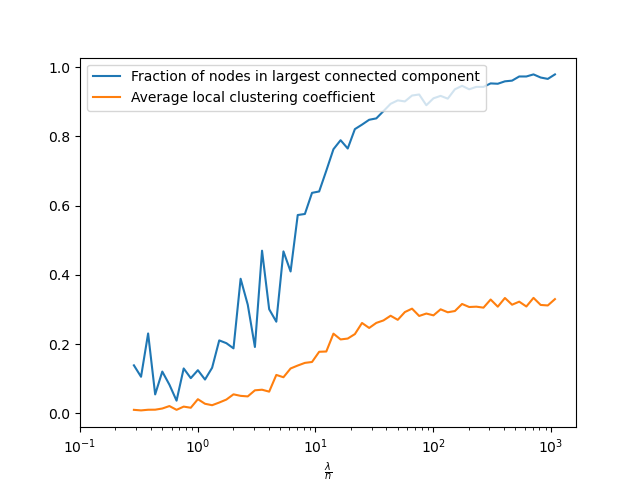}
		\caption{Largest component and  clustering}
		\label{fig:alpha_beta_local_clustering}
	\end{subfigure}
	\caption{Clustering versus degree and the largest component size}
	\label{fig:clustering degree component}
\end{figure}

In Figure \ref{fig:clustering degree component}  (a) we examine several clustering characteristics of  the stationary network gene\-ra\-ted by the ``general triadic model'' with $\mu=15000$ and $\lambda=20000$. 
Given integer $k\ge 2$, let $g(k)$ denote the number of vertices $v$ of degree $d(v)=k$. Let $f(k)=\frac{1}{g(k)}\sum_{v:\, d(v)=k}C^{\text{\tiny L}}_v(G)$ denote the average value of the local clustering coefficient over the set of vertices of degree $k$. We put $f(k)=0$ for $g(k)=0$. We call $f$ the ``local clustering coefficient curve''. The fact that $f$ is decreasing tells us that the local clustering coefficient  negatively correlates with  vertex degree, a phenomenon observed in many sparse real networks (\cite{Foudalis_J_P_Sideri_2011},
\cite{Ravasz_S_M_O_Barabasi_2002},
\cite{Ravasz_Barabasi_2003},  
\cite{Vazquez_Pastor-Satorras_Vespignani_2002}). The ``general triadic model'' reproduces this network property. We also mention that the edge density $0.004$ is by two orders less than the average local clustering coefficient. Hence the network is sparse and highly clustered.

Lastly, we  touch on the question of the component structure. One may wonder whether the high values of the clustering coefficients are caused by  a  few (perhaps one) relatively small, but dense subgraphs. Figure 
\ref{fig:clustering degree component} (b) shows that this is not the case. The stationary network generated by ``general triadic model'' admits a large connected component collecting a fraction of nodes.  For simplicity we put $\mu\equiv 0$ (no triad protection).  Hence the only remaining  parameter to vary is $\lambda$.
On the horizontal axis we depict values of $\frac{\lambda}{n}$. We recall that the number of vertices $n=1000$ remains fixed.  

 \subsection{Rigorous results}

Let $f$ be a real valued function defined on the set of graphs with 
the vertex set $V$. For example, it can be the number of edges  
$f(G)=|E|$ of graph $G=(V,E)$, or  the number of  triangles 
$f(G)=N_{\Delta}(G)$, etc.
 For a stationary Markov chain $\bf G$ the function   
$t\to\E f(G_t)$ is a constant. Hence $\frac{\partial}{\partial t}\E f(G_t)=0$.
This identity, when applied to properly chosen function $f$,  can 
give useful information about  average characteristics of the network at stationarity.
 We explore two instances. Choosing $f(G)=|E|$ we show lower and upper bounds for  the average edge density $e_t:=\E e(G_t)$; choosing $f(G)=N_{\Delta}(G)$ we infer about 
the number of triangles. 

Since for stationary ${\bf G}$
the average edge density
 $e_t$ and average clustering coefficient    $\E {\bar C}^{\text{\tiny L}}(G_t)$
do not depend on $t$, we sometimes drop the subscript $t$ and write $e=e_t$ and  ${\bar C}^{\text{\tiny L}}=\E {\bar C}^{\text{\tiny L}}(G_t)$.
We observe that, by symmetry, the probability distribution of  bivariate random variable 
$(d_v(G_t),\Delta_v(G_t))$ is the same for all $v\in V$. Furthermore,  for a stationary network this distribution does not depend   on $t$ either. We denote by $(d,\Delta)$ a bivariate random variable having the same distribution as 
$(d_v(G_t),\Delta_v(G_t))$.

To make calculations feasible we assume for the rest of the section that  the products
$\lambda_i\lambda_j$ and
$\mu_i\mu_j$ in (\ref{model1}) do not depend on $i,j$. 
In this case (\ref{model1}) reads as follows 
\begin{align}\label{modelUW}
a_{ij}(G) =
\begin{cases}
\lambda_0+\lambda\nu_{ij}(G, \alpha)
\qquad
\quad
\
\
 {\text{for}}
\quad 
\{i,j\}\not\in E,
\\
\left(\mu_0-\mu\nu_{ij}(G,\beta)\right)_+
\qquad
{\text{for}}
\quad
\{i,j\}\in E,
\end{cases}
\end{align}
where 
$\lambda,\lambda_0>0$ and $\mu,\mu_0>0$. Below ${\bf G}$ denotes 
a stationary Markov chain defined by (\ref{modelUW}).

{\it Edge density.} 
We have  that
\begin{align}
\label{2024-08-05}
e\ge \frac{\lambda_0}{\lambda_0+\mu_0}.
\end{align}
For $\alpha,\beta\ge 2$ we have
 \begin{align}
\label{2024-09-23}
 e\le \frac{\lambda_0+\frac{1}{n-1}\max\{\lambda,\mu\}}{\lambda_0+\mu_0}.
\end{align}
For $\alpha,\beta\ge 1$  and  
$\lambda_0+\mu_0>\max\{\lambda,\mu\}$ we have that
\begin{align}
\label{2024-08-06}
 e\le \frac{\lambda_0}{\lambda_0+\mu_0-\max\{\lambda,\mu\}}.
\end{align}

An important conlcussion  to draw from inequalities  (\ref{2024-08-05}), (\ref{2024-09-23}), 
(\ref{2024-08-06}) 
is that for $\frac{\mu_0}{\lambda_0}$ of the order $n$ and 
$\max\{\lambda,\mu\}$ of the order $\mu_0$   the network $G_t$  
is sparse and has average edge 
density of the order $n^{-1}$ as $n\to+\infty$.

  
Proof of (\ref{2024-08-05}), (\ref{2024-09-23}), and (\ref{2024-08-06}). 
Equation $\frac{\partial}{\partial t}\E |E_t|=0$
implies 
\begin{align}\label{2024-08-22}
\E\sum_{\{i,j\}\notin E_t}a_{ij}(G_t)
=
\E\sum_{\{i,j\}\in E_t}a_{ij}(G_t).
\end{align}
In view of (\ref{modelUW}) we can write the latter identity in the form
\begin{align}
\label{2024-08-05_1}
\E\left(
\lambda_0\left({{n}\choose{2}}-|E_t|\right)
+
\nu'_t +\nu''_t
-
\mu_0
|E_t|
\right)
=0.
\end{align}
where 
\begin{align}
\nu'_t
=
\lambda
\sum_{\{i,j\}\notin E_t}
\sum_{v\in N_{ij}}\frac{1}{d_v^{\alpha}}
\qquad
{\text{and}}
\qquad
\nu''_t
=
\sum_{\{i,j\}\in E_t}
\min\left\{\mu_0,\,
\mu \sum_{v\in N_{ij}}
\frac{1}{d_v^{\beta}}
\right\}
\end{align}
account for the contribution of the clustering weights 
$\lambda\nu_{ij}(G_t,\alpha)$ and $\mu\nu_{ij}(G_t,\beta)$. 
Here we write,
for short, $N_{ij}=N_{ij}(G_t)$  and $d_w=d_w(G_t)$. 
By the linearity of expectation, we obtain from (\ref{2024-08-05_1}) that 
\begin{align}
\label{2024-08-05_2}
\lambda_0-(\lambda_0+\mu_0)e_t
+
{{n}\choose{2}}^{-1}
\E(\nu'_t+\nu''_t)
=0.
\end{align}
 The inequalities $\nu'_t\ge 0$, $\nu''_t\ge 0$ 
imply $\lambda_0-(\lambda_0+\mu_0)e_t\ge 0$. We arrived to
 lower bound 
 (\ref{2024-08-05}).

Let us show upper bounds (\ref{2024-09-23}), (\ref{2024-08-06}). We denote $\tau:=\min\{\alpha,\beta\}$ and 
estimate 
\begin{align*}
\nu'_t
&
\le
\lambda
\sum_{\{i,j\}\notin E_t}
\sum_{v\in N_{ij}}\frac{1}{d^{\tau}_v},
\\
\nu''_t
&
\le
\mu 
\sum_{\{i,j\}\in E_t}
\sum_{v\in N_{ij}}
\frac{1}{d_v^{\beta}}
\le
\mu 
\sum_{\{i,j\}\in E_t}
\sum_{v\in N_{ij}}
\frac{1}{d^{\tau}_v}.
\end{align*}
Combining these inequalities we obtain
\begin{align*}
\nu'_t+\nu''_t
&\le 
\max\{\lambda,\mu\}
\sum_{1\le i<j\le n}
\sum_{v\in N_{ij}}\frac{1}{d_v^\tau}
=
\max\{\lambda,\mu\}
\sum_{v\in V: \, d_v\ge 2}
\frac{1}{d^{\tau}_v}{{d_v}\choose{2}}
\\
&
\le
 \max\{\lambda,\mu\}
 \frac{1}{2}\sum_{v\in V}d^{2-\tau}_v.
\end{align*}
For $\tau=2$ we have 
$\nu'_t+\nu''_t\le \frac{n}{2} \max\{\lambda,\mu\}$. Invoking this inequality in (\ref{2024-08-05_2}) we obtain 
(\ref{2024-09-23}).
For $\tau=1$ we have $\nu'_t+\nu''_t\le \frac{1}{2}\max\{\lambda,\mu\}|E_t|$.
Now (\ref{2024-08-05_2}) yields 
(\ref{2024-08-06}).

{\it Special case of $\alpha=2$.} In this special case we consider a slightly 
modified version of (\ref{modelUW}) that includes the ``correction term''
\begin{displaymath}
\varkappa_{ij}(G)
=
\frac{1}{n-1}
\left(
{\mathbb I}_{\{d_i(G)=0\}}
+
{\mathbb I}_{\{d_j(G)=0\}}
\right)
+
\frac{1}{n-2}
\left(
{\mathbb I}_{\{d_i(G)=1\}}
+
{\mathbb I}_{\{d_j(G)=1\}}
\right).
\end{displaymath} 
In addition, we replace  $\nu_{ij}(G,2)$ by related quantity $\nu^*_{ij}(G)=\sum_{w\in N_{ij}(G)}{\binom{d_w(G)}{2}}^{-1}$.
We  set
\begin{align}\label{Dominykovariantas}
a_{ij}(G) =
\begin{cases}
\lambda_0+\lambda\nu^*_{ij}(G)+\lambda\varkappa_{ij}(G)
\qquad
\
\
 {\text{for}}
\quad 
\{i,j\}\not\in E,
\\
\mu_0
\qquad
\qquad
\qquad
\qquad
\qquad
\quad
\
\
\
{\text{for}}
\quad
\{i,j\}\in E.
\end{cases}
\end{align}
The reason for such  a modification is that it admits  a closed form expression for the average edge density.

For a stationary network Markov chain defined by (\ref{Dominykovariantas}) we have that
\begin{align}\label{2024-087-21+3}
e=\frac{\lambda_0+\frac{2}{n-1}\lambda(1-{\bar C}^{\text{\tiny L}})}{\lambda_0+\mu_0}
\end{align}
Noting that ${\bar C}^{{\text{\tiny L}}}\le 1$  we obtain from 
(\ref{2024-087-21+3})  the upper and lower bounds for the average edge density
\begin{displaymath}
\frac{\lambda_0}{\lambda_0+\mu_0}
\le 
e
\le 
\frac{\lambda_0+\frac{2}{n-1}\lambda}{\lambda_0+\mu_0}.
\end{displaymath}
Letting $n\to+\infty$ and choosing 
  $\frac{\mu_0}{\lambda_0}$, $\frac{\lambda}{\lambda_0}$ 
  and $\frac{\mu}{\lambda_0}$ of the order $n$ we have  that  $e$ is of the order $n^{-1}$. Hence the model produces a sparse dynamic network. 

Proof of (\ref{2024-087-21+3}). 
Equation (\ref{2024-08-22}) 
implies 
\begin{align}
\label{2024-08-22+1}
\E 
\sum_{\{i,j\}\notin E_t}
\left(
\lambda_0
+
\lambda\varkappa_{ij}
+
\lambda \sum_{w\in N_{ij}}\frac{1}{\binom{d_w}{2}}\right)
=\mu_0\E|E_t|.
\end{align}
Here we write, for short, $\varkappa_{ij}=\varkappa_{ij}(G_t)$, $N_{ij}=N_{ij}(G_t)$ and $d_w=d_w(G_t)$. Invoking the identities
\begin{align*}
&
\sum_{\{i,j\}\notin E_t}1
=\binom{n}{2}-|E_t|,
\\
&
\sum_{\{i,j\}\notin E_t}\varkappa_{ij}
=
\sum_{w\in V}{\mathbb I}_{\{d_w= 0\}}
+
\sum_{w\in V}{\mathbb I}_{\{d_w=1\}},
\\
&
\sum_{\{i,j\}\not\in E_t}\sum_{w\in N_{ij}}\frac{1}{\binom{d_w}{2}}
=
\sum_{w\in V: \, d_w\ge 2}
\frac{\binom{d_w}{2}-\Delta_w(G_t)}{\binom{d_w}{2}}
=
\sum_{w\in V: \, d_w\ge 2}(1-C^{\text{\tiny L}}_w(G_t))
\end{align*}
and dividing both sides of (\ref{2024-08-22+1}) by $\binom{n}{2}$ we have
\begin{displaymath}
\lambda_0(1-e)
+
\frac{2\lambda}{n-1} 
\left(
\PP\{d\le 1\}+\PP\{d\ge 2\}-
{\bar C}^{\text{\tiny L}}
\right)=\mu_0 e,
\end{displaymath}
where $d$ denotes the degree of a randomly selected vertex.
We have arrived to (\ref{2024-087-21+3}).

{\it Number of triangles.}  
For a stationary network Markov chain defined by (\ref{modelUW}) where $0<\alpha\le 2$
we show that
\begin{align}\label{2024-08-21}
\E \Delta
\ge 
\frac{\lambda}{4(\lambda_0+\mu_0+\lambda)}\PP\{d\ge 2\}.
\end{align}
An important conlcussion  to draw from inequality
 (\ref{2024-08-21})  
 is   that choosing 
  $\frac{\mu_0}{\lambda_0}$, $\frac{\lambda}{\lambda_0}$ 
  and $\frac{\mu}{\lambda_0}$ of the order $n$ one  can obtain a sparse
  stationary dynamic network with the property that the average number of triangles incident to a vertex of degree at least two (formally, the conditional expectation $\E (\Delta|d\ge 2)=\frac{\E \Delta}{\PP\{d\ge 2\}}$)
   is bounded from below by a constant. Note that $d_v(G_t)\le 1$ implies $\Delta_v(G_t)=0$. Hence $\Delta_v(G_t)=\Delta_v(G_t){\mathbb I}_{\{d_v(G_t)\ge 2\}}$ and $\Delta=\Delta{\mathbb I}_{\{d\ge 2\}}$. 

Proof of (\ref{2024-08-21}).
Equation
$\frac{\partial}{\partial t}N_{\Delta}(G_t)=0$
implies
\begin{align}
\label{2024-08-07}
\E \sum_{\{i,j\}\notin E_t}|N_{ij}(G_t)|a_{ij}(G_t)
=
\E \sum_{\{i,j\}\in E_t}|N_{ij}(G_t)|a_{ij}(G_t).
\end{align}
Here the left sum evaluates the average birth rate of triangles: connecting a pair of non-adjacent vertices  $i,j$ by an edge creates  
$|N_{ij}(G_t)|$ new triangles. The right sum  evaluates the average death rate of triangles: deletion of an edge $\{i,j\}\in E_t$ eliminates 
$|N_{ij}(G_t)|$ triangles from $G_t$.
Furthermore, for $\{i,j\}\in E_t$ we have 
$a_{ij}(G_t)\le \mu_0$. Hence  the sum on the right of 
(\ref{2024-08-07})
\begin{align}
\label{2024-08-21+1}
\sum_{\{i,j\}\in E_t}|N_{ij}(G_t)|a_{ij}(G_t)
\le
\sum_{\{i,j\}\in E_t}|N_{ij}(G_t)|\mu_0
=
\sum_{v\in V:\, d_v(G_t)\ge 2}\Delta_v(G_t)\mu_0.
\end{align}
In the last identity we use the observation that 
$\Delta_v(G_t)$ counts 
edges whose both endpoints are adjacent to $v$. Similarly for the sum on the left of 
(\ref{2024-08-07})
\begin{align}
\label{2024-08-21+2}
\sum_{\{i,j\}\notin E_t}|N_{ij}(G_t)|a_{ij}(G_t)
\ge
 \sum_{v\in V: \, d_v(G_t)\ge 2}
\left(
{{d_v(G_t)}\choose{2}}-\Delta_v(G_t)
\right)\left(\lambda_0+\frac{\lambda}{d_v^{\alpha}(G_t)}
\right)
\end{align}
Here we  use the observation that ${{d_v(G_t)}\choose{2}}-\Delta_v(G_t)$ counts pairs 
$\{i,j\}$ of neighbours of $v$ that are non-adjacent ($\{i,j\}\notin E_t$). Inequality (\ref{2024-08-21+2})  
follows from the fact that
$a_{ij}(G_t)\ge \lambda_0+\frac{\lambda}{d_v^{\alpha}}$
 for  each $v\in N_{ij}(G_t)$. 

Invoking (\ref{2024-08-21+1}) and (\ref{2024-08-21+2}) in (\ref{2024-08-07}) we obtain
\begin{displaymath}
\E
\sum_{v\in V: \, d_v(G_t)\ge 2}
\left({{d_v(G_t)}\choose{2}}-\Delta_v(G_t)\right)
\left(\lambda_0+\frac{\lambda}{d_v^{\alpha}(G_t)}\right)
\le
\E
\sum_{v\in V:\, d_v(G_t)\ge 2}
\Delta_v(G_t)\mu_0.
\end{displaymath}
Recall that the probability distribution of bivariate random variable 
$(d_v(G_t), D_v(G_t))$ is the same for all $v\in V$. Collecting the  terms $\Delta_v(G_t)$ on the right  and dividing both sides by $n$ we have
\begin{displaymath}
\frac{\lambda_0}{2}\E\bigl(d(d-1)\bigr) 
+
\frac{\lambda}{2}\E (d^{1-\alpha}(d-1){\mathbb I}_{\{d\ge 2\}}
\le
(\lambda_0+\mu_0)\E \Delta
+
\lambda \E 
\left(\Delta_v 
\frac{{\mathbb I}_{\{d_v\ge 2\}}}{d_v^{\alpha}}\right).
\end{displaymath}
Next we upper bound 
$\lambda \E \left(\Delta_vd_v^{-\alpha}
{\mathbb I}_{\{d_v\ge 2\}}\right)
\le 
\lambda \E \Delta $ 
and obtain
\begin{displaymath}
\frac{\lambda}{2}
\E 
\left(
d^{1-\alpha}(d-1){\mathbb I}_{\{d\ge 2\}}
\right)
\le
(\lambda_0+\mu_0+\lambda)\E \Delta.
\end{displaymath}
Furthermore, using inequality
$\frac{1}{2}\le \frac{d-1}{d} \le \frac{d-1}{ d^{\alpha-1}}$, which holds  for $0<\alpha\le 2$ and $d\ge 2$ we lower bound the left side  by 
$\frac{\lambda}{4}$ and obtain  inequality equivalent to (\ref{2024-08-21})
\begin{displaymath}
\frac{\lambda}{4}\PP\{d\ge 2\}
\le 
(\lambda_0+\mu_0+\lambda)\E \Delta.
\end{displaymath}

\section{Dynamic affiliation network}

Let
${\bf H}=\{H_t=(V\cup W ,E_t),\, t\ge 0\}$ be a continuous time Markov 
chain, whose states are bipartite graphs with the bipartition  $V\cup W$,
where  $V=\{1,\dots, n\}$ and $W=\{1,\dots, m\}$. 
Transitions of ${\bf H}$ are defined as follows. 
Given $H_t$ (the state occupied at time $t$), the update takes place  at time  $t':=t+\min_{(i,u)\in V\times W} B_{iu}$
when  the pair $(i_0,u_0)$ with 
$B_{i_0u_0}=\min_{(i,u)\in V\times W} B_{iu}$  changes its adjacency status.
Here $B_{iu}=B_{iu}(H_t)$, 
are independent exponential waiting times with intensities $b_{iu}=b_{iu}(H_t)$ defined below.   Thus, at time $t'$  chain ${\bf H}$ 
jumps to the state $H_{t'}=(V\cup W, E_{t'})$, where the edge sets $E_t$ and $E_{t'}$ differ in the single edge $i_0\sim u_0$.
Markov chain ${\bf H}$ defines dynamic affiliation network 
${\bf G}'=\{G'_t =(E'_t,V), t\ge 0\}$:
 for each $t$  any two vertices $i,j\in V$  are adjacent in $G'_t$ whenever $i$ and $j$ 
have  a common neighbour 
in $H_t$.

Now we  define intensities $b_{iu}$. We fix $\mu>0$ and  assign positive weights $y_i$ and $x_u$  to
$i\in V$ and $u\in W$ that model activity of actors and attrativeness of attributes.  For a bipartite graph $H=(V\cup W, E)$  we  set
\begin{align}\label{model}
b_{iu}(H) =
\begin{cases}
y_ix_u
\qquad
 {\text{for}}
\quad 
(i,u)\not\in E,
\\
\mu
\quad
\qquad
\,
{\text{for}}
\quad
(i,u)\in E.
\end{cases}
\end{align}
 Clearly, ${\bf H}$ has a
 unique stationary distribution defined by the weight 
 sequences $\{y_i\}_{i=1}^n$, $\{x_u\}_{u=1}^m$ and $\mu$.
Furthermore,  $\bf H$ comprises of 
$n\times m$ independent continuous Markov chains describing  adjacency  
dynamic 
of
each vertex pair  $(i,u)\subset V\times W$ separately, where 
the Markov chain of a pair
 $(i,j)$ has two states $i\sim u$ and $i\not\sim u$
 whose
 exponential holding times have intensities $\mu$ and $y_ix_u$ respectively.  Thus, at stationarity, a snaphot $H_t$ represents a random bipartite graph, where edges are inserted independently  with probabilities  
 \begin{align}
 \label{2024-08-27}
 \PP\{i\sim u\}= \frac{y_ix_u}{y_ix_u+\mu}=:p_{iu},
 \end{align}
 for $(i,u)\in V\times W$. 
We assume in what follows that dynamic affiliation network ${\bf G}'$
is 
 defined by a stationary Markov chain ${\bf H}$ satisfying  
 (\ref{2024-08-27}). In this case probability distributions of random graphs $G'_t$ and $H'_t$ do not depend on $t$ and with a little abuse of notation we write, for short, $G'=G'_t$ and $H=H_t$. 
 %
 We show  that 
 $G'$ admits tunable degree distribution and 
nonvanishing global clustering coefficient. 

We will use the following notation. 
By 
$P_{y,n}=\sum_{i=1}^n\delta_{y_i}$ 
and 
$P_{x,m}=\sum_{u=1}^m\delta_{x_u}$ 
we denote empirical distributions of the sequences 
$\{y_1,\dots, y_n\}$ and
$\{x_1,\dots, x_m\}$. Here $\delta_t$ stands for the degenerate distribution that assigns mass $1$ to point $t$.
Furthermore we denote 
\begin{displaymath}
\langle x^s\rangle=\frac{1}{m}\sum_{u\in [m]}x_u^s,
\quad
\langle y^s\rangle=\frac{1}{n}\sum_{i\in [n]}y_i^s,
\quad
\gamma^2=\frac{m}{n},
\quad
\varkappa=\frac{nm}{\mu^2}.
\end{displaymath}
It is important to mention that  the ratio 
$\gamma^2=\frac{m}{n}$ correlates negatively with the clustering  strength. More precisely,
the global clustering coefficient of $G'$
is asymptotically inversely proportional
to  $\gamma$ for large  $n$ and $m$,   see (\ref{2024-10-29}) below.

{\it Degrees of $G'$}. Here we show that
the expected value $\E d_i$ of the degree $d_i=d_i(G')$ of vertex $i$ is  
approximately proportional to  its weight
$y_i$.  Moreover,   $d_i$  has asymptotic compound Poisson distribution as the network size $n\to+\infty$.

\begin{tm}\label{T1}
For  each $i\in V$ we have 
\begin{align}
\label{2024-09-07}
0
\le
 y_i \varkappa\,\langle x^2\rangle\langle y\rangle - \E d_i
 \le 
\frac{\varkappa}{\mu}y_i
 \left(\langle  x^3\rangle\langle y^2\rangle
 +
 y_i\langle x^3\rangle
 \langle y\rangle
 \right)
 +
\frac{\varkappa^2}{n} 
y_i^2\langle x^2\rangle^2\langle y^2\rangle
 +
\frac{1}{n}y_i\langle x^2\rangle.
\end{align}
\end{tm}

It follows from Theorem \ref{T1} that for large $n,m$  and 
$\mu=\mu(n,m)$ of the order $\sqrt{nm}$ the expected degree of 
a vertex $i$ in $G'$  
is asymptotically proportional its ``activity'' weight $y_i$. 

\begin{cor}\label{C1}
Let  $n,m\to+\infty$. Put $\mu=\sqrt{nm}$ and  assume that for some $c>0$
we have 
$\langle x^3\rangle\le c$  
and $\langle y^2\rangle\le c$  uniformly in $n,m$. Then  for each  $i$  
\begin{align}
\label{2024-09-07+1}
\E d_i
=
y_i\langle x^2\rangle\langle y\rangle+O\left(\frac{y_i^2}{\sqrt{nm}}\right).
\end{align}
\end{cor}

In the sparse regime, when $\E d_i$ remains bounded as $n,m\to+\infty$, the probability distribution of $d_i$ does not concentrate around the expected value $\E d_i$.  Theorem \ref{T2} below shows  that $d_i$ has a compound Poisson asymptotic distribution.
Recall that compound Poisson distribution is the probability 
distribution of a randomly stopped sum $\sum_{k=1}^{\Lambda}\xi_k$, 
where $\xi_1,\xi_2,\dots$ are independent and identically distributed 
random variables, which are independent of  Poisson  random 
variable $\Lambda$. We write 
$\Lambda\sim {\cal P}(\lambda)$, where $\lambda:=\E \Lambda$ denotes the 
expected value and denote by ${\cal CP}(\lambda,P_{\xi})$ the (compound Poisson) distribution of 
$\sum_{k=1}^{\Lambda}\xi_k$. Here $P_{\xi}$ denotes the (common) 
probability distribution of $\xi_k$.

Let $x_1,x_2,\dots$ and $y_1,y_2,\dots$ be positive infinite sequences of weights. In Theorem \ref{T2} we consider random affiliation networks  $G'_{n,m}$, $n,m=1,2,\dots$, based on  respective bipartite random graphs $H_{n,m}$ whose edges are inserted independently with probabilities  (\ref{2024-08-27}). Note that each $H_{n,m}$  is defined by 
truncated (finite)  sequences $\{x_1,\dots, x_m\}$ and $\{y_1,\dots, y_n\}$.

To formulate our next result we need the following conditions: for $n,m\to+\infty$ we have

(i) $P_{x,m}$ converges weakly to some probability distribution, say $P_X$,  having  a finite first moment
$\int sP_X(ds)<\infty$ and $\langle x\rangle$ converges to $\int sP_X(ds)$;

(ii) the family of distributions $\{P_{y,n},\, n=1,2,\dots\}$ is uniformly integrable and $\langle y\rangle$ converges to some number $a_y>0$.

\begin{tm}\label{T2} Let $\mu=\sqrt{nm}$. Let $n\to+\infty$. Assume that $m=m(n)$ is such that  $m/n$ converges to some $\gamma_o>0$. Assume that (i) and (ii) hold. Denote 
$a_x= \int sP_X(ds)$ and introduce function 
$s\to \lambda_s=sa_y\gamma_o^{-1}$. 
For each $i=1,2,\dots$ the probability distribution of $d_i$ converges weakly to 
the compound Poisson distribution ${\cal CP}(y_ia_x\gamma_o,Q)$, where 
the discrete probability distribution $Q$  assigns probabilities 
\begin{align}
\label{T2+1}
Q(\{t\})
=
\int 
\frac{s}{a_x}
e^{-\lambda_s}\frac{\lambda_s^t}{t!} P_X(ds).
\end{align}
to integers $t=0,1,2,\dots$. 
\end{tm}
We note that $Q$ is a mixture of Poisson distributions. To sample from 
$Q$ one can use the two step procedure: 1)  generate a (size biased) random variable  ${\tilde X}$ according 
to the distribution $\PP\{{\tilde X}=s\}=\frac{s}{a_x}\PP\{X=s\}$, 
$s=0,1,\dots$;
2) sample Poisson random variable with rate
${\tilde X}a_y\gamma_o^{-1}$.

\bigskip

{\it Clustering in $G'$}.  We recall that $N_{\Delta}(G')$ denotes the number of 
triangles in $G'$ and $N_{\Lambda}(G')$ denotes the number of $2$-
paths in $G'$. 

\begin{tm}\label{T3+} Let $\mu=\sqrt{nm}$. Let $n\to+\infty$. Assume that $m=m(n)$ is such that  $m/n$ converges to some $\gamma_o>0$.
Assume that for some constant 
$c>0$  we have  $\langle x^5\rangle <c $ and $\langle y^4\rangle <c$ 
for all $n$. Then
\begin{align}
\label{2024-10-15+1}
N_{\Delta}(G')
&
=
\frac{n}{6\gamma}\langle x^3\rangle \langle y\rangle^3+
o_P(\sqrt{n}),
\\
\label{2024-10-26+30}
N_{\Lambda}(G')
&
=
\frac{n}{2\gamma}\langle x^3\rangle \langle y\rangle^3
+
\frac{n}{2}
\langle x^2\rangle^2
\langle y^2\rangle
\langle y\rangle^2
+
O_P(\sqrt{n}).
\end{align}
In particular, the global clustering coefficient 
\begin{align}
\label{2024-10-29}
C^{\text{\tiny GL}}(G)
=
\frac
{\langle x^3\rangle \langle y\rangle^3}
{\langle x^3\rangle \langle y\rangle^3
+
\gamma_o\langle x^2\rangle^2 \langle y^2\rangle\langle y\rangle^2
}
+
o_P(1).
\end{align}
\end{tm}
We remark that  conditions $\langle x^5\rangle <c$ and $\langle y^4\rangle <c$ of Theorem \ref{T3+}
can be relaxed. We expect that the minimal conditions  $\langle x^3\rangle <c$ and $\langle y^2\rangle <c$ plus  the uniform integrability 
of  $t^3P_{x,m}(dt)$ and 
$t^2P_{y,n}(dt)$ would suffice.

\section {Concluding remarks} 
We presented two dynamic network models that  generate sparse and clustered stationary  networks. Both models seems natural as they mimic dynamics of real network processes. 
Luckily, for  rigorous analysis of dynamic  affiliation network we can use
 techniques  developed for random intersection graphs
\cite {B_G_J_K_R_2015_models},
\cite{Grohn_Karjalainen_Leskela_2024}, \cite{Hofstad_Komjathy_Vadon_2021},
\cite{Hofstad_book_2_2024},
\cite{Kurauskas2022},
\cite{Milewska_Hofstad_Zwart}. 
 On the other hand  we have only a  few rigorous results for
 stationary  Markov chains with clustering like (\ref{model1}), (\ref{modelUW}). 
 %
%
 It would be interesting to learn more about network structure and properties
 of this model via rigorous analysis. 

\medskip

{\bf Acknowledgement}. 
The diagrams  were generated using Matplotlib \cite{Maplotlib}. 
Authors thank Information technology research center of Vilnius University for a high performance computing resourses.

\section{Appendix}

Here we prove Theorems \ref{T1}, \ref{T2}, \ref{T3+}. Before proofs we introduce some notation.

For $v\in V\cup W$ we denote by $N_v$ the set of neighbours of $v$ in $H$. Note that for $i\in V$ and $u\in W$ we have   $N_i\subset W$ and 
$N_u\subset V$.
For $i\in V$, $u\in W$ we denote by ${\mathbb I}_{iu}={\mathbb I}_{\{i\sim u\}}$
 the indicator function of the event $i\sim u$ 
 (meaning that $i$ and $u$ are adjacent in $H$).
For $i\in V$ we denote by $N_i^{G}$ the set of neighbours of $i$ in $G'$;
$d_i=|N_i^{G}|$ denotes the degree of $i$ in $G'$.

We write, for short, 
\begin{align*}
S_a(x)=\sum_{i\in[m]}x_i^a=m\langle x^a\rangle,
\qquad
S_a(y)=\sum_{j\in[n]}y_j^a=n\langle y^a\rangle,
\qquad
p_{ik}^{\star}=\frac{y_ix_k}{\mu}.
\end{align*}

In the proof (sometimes without mentioning) we  apply   inequalities
 \begin{align}
 \label{2024-09-10}
p_{iu}^{\star}
\ge  
p_{iu}
\ge 
\left(1-p_{iu}^{\star}\right)
p_{iu}^{\star},
 \quad
\prod_{j}p_{i_ju_j}^{\star}
\ge
\prod_{j}p_{i_ju_j}\ge \left(1-\sum_jp_{i_ju_j}^{\star}\right)
\prod_{j}p_{i_ju_j}^{\star}.
\end{align}
The first  (and third) inequality is obvious. The second one follows from the inequalities  $\frac{a}{b}\ge \frac{a}{a+b}\ge \left(1-\frac{a}{b}\right)\frac{a}{b}$,
 for $a>0$ and $b>0$.  The fourth one is obtained by iterating the second inequality.

{\it Proof of Theorem \ref{T1}.}
We only prove (\ref{2024-09-07}).
Let
 \begin{align*}
 d_i^{(1)}
&
 =
  \sum_{u\in N_i}|N_u\setminus\{i\}|
  =
  \sum_{u\in W}{\mathbb I}_{iu}
 \sum_{j\in [n]\setminus\{i\}}{\mathbb I}_{ju},
 \\
 R_i
 &
 =
 \sum_{\{u,v\}\subset N_i}\bigl|(N_u \setminus\{i\})\cap( N_v\setminus\{i\})\bigr|
 =
 \sum_{\{u,v\}\subset W}
 {\mathbb I}_{iu}
 {\mathbb I}_{iv}
 \sum_{j\in [n]\setminus \{i\}}
  {\mathbb I}_{ju}
 {\mathbb I}_{jv}.
 \end{align*}
 Using inclusion-exclusion inequalities we estimate the number of elements of the set 
 $N_i^{G}=\cup_{u\in N_i}(N_u\setminus\{i\})$. 
 We have
 \begin{align}
 \label{2024-09-30}
 d_i^{(1)}-R_i
 \le 
 |N_i^{G}|
 \le 
 d_i^{(1)}.
 \end{align}
 To prove the left inequality of
 (\ref{2024-09-07}) 
 we use $\E d_i\le \E d_i^{(1)}$ (see (\ref{2024-09-30})) and invoke
 bound (\ref{2024-09-30+1}) shown below.
 To show the right inequality of (\ref{2024-09-07}) 
 we use $\E d_i^{(1)}-\E R_i\le \E d_i$ (see (\ref{2024-09-30})) and invoke  
 bounds (\ref{2024-09-30+2}), (\ref{2024-09-30+3}) below.
 Finally, using (\ref{2024-09-10}) we estimate $\E d_i^{(1)}$ and $\E R_i$.
\begin{align}
\label{2024-09-30+1}
\E d_i^{(1)}
&
=\sum_{u\in [m]}p_{iu}\sum_{j\in[n]\setminus\{i\}}p_{ju}
\le
\sum_{u\in [m]}p_{iu}\sum_{j\in[n]}p_{ju}
\le
\frac{mn}{\mu^2}
 y_i
 \langle x^2\rangle
 \langle y\rangle,
 \\
\nonumber
 \E d_i^{(1)}
 &
 =
 \sum_{u\in [m]}\sum_{j\in[n]}p_{iu}p_{ju}-\sum_{u\in [m]}p_{iu}^2
\\
\nonumber
&
 \ge
 \sum_{u\in [m]}\sum_{j\in[n]}
 p_{iu}^{\star}p_{ju}^{\star}
 \left(1-p_{iu}^{\star}-p_{ju}^{\star}\right)
-
\sum_{u\in [m]}(p_{iu}^{\star})^2
\\
 \label{2024-09-30+2}
&
=
\frac{mn}{\mu^2}
y_i
\langle x^2\rangle
\langle y\rangle
-
\frac{mn}{\mu^3}
y_i^2
\langle x^3\rangle
\langle y\rangle
-
\frac{mn}{\mu^3}
y_i
\langle x^3\rangle
\langle y^2\rangle
-
\frac{m}{\mu^2}
y_i^2
\langle x^2\rangle
\end{align}
and
\begin{align}
\nonumber
\E R_i
&
\le
\sum_{\{u,v\}\subset [m]}p_{iu}p_{iv}\sum_{j\in [n]}p_{ju}p_{jv}
\le
\sum_{\{u,v\}\subset [m]}p_{iu}^{\star}p_{iv}^{\star}
\sum_{j\in [n]}p_{ju}^{\star}p_{jv}^{\star}
\\
 \label{2024-09-30+3}
&
\le\frac{1}{2} \frac{m^2n}{\mu^4}
y_i^2
\langle x^2\rangle^2
\langle y^2\rangle.
\end{align}
Proof of Theorem \ref{T1} is complete.

\bigskip

{\it Proof of Theorem \ref{T2}.} 
Before the proof we introduce some notation and collect auxiliary results.
Given two random variables $\xi$ and $\zeta$ we denote by 
$d_{TV}(P_\xi,P_\zeta)$ the total variation distance between the probability 
distributions $P_\xi$ and $P_\zeta$ of $\xi$ and $\zeta$. With a little abuse of notation we also write $d_{TV}(\xi,\zeta)$.

We fix vertex  $i\in V$.  For $u\in W$ we write, for short,  
$\zeta_u={\mathbb I}_{iu}$ and
$\xi_u=\sum_{j\in V\setminus \{i\}}{\mathbb I}_{ju}$.
Let $\zeta_u^*$, $\zeta_u^{\star}$, $\zeta_u^{o}$ and
$\xi_u^*$, $\xi_u^{\star}$, $\xi_u^{o}$ be Poisson random variables with expected values 
\begin{align*}
&
\E \zeta^*_u=p_{iu},
\quad
\
\E \zeta^{\star}_u=p_{iu}^{\star},
\quad
\
\E \zeta^{o}_u=\frac{\gamma_o}{\gamma}p_{iu}^{\star},
\\
&
\E \xi^*_u=\sum_{j\in V\setminus\{i\}}p_{ju},
\quad
\
\E \xi^{\star}_u=\sum_{j\in V\setminus\{i\}}p_{ju}^{\star},
\quad
\
\E \xi^{o}_u=\gamma_o^{-1}x_ua_y.
\end{align*}
Note that $\E \zeta^*_u=\E \zeta_u$ and $\E \xi^*_u=\E \xi_u$.
Let 
$(\xi_u(k),\xi_u^*(k),\xi_u^{\star}(k),\xi_u^{o}(k)), k\ge 1$, be iid copies 
of $(\xi_u,\xi_u^*,\xi_u^{\star}, \xi_u^{o})$. We assume that each   collection 
\begin{align*}
&
\{\zeta_u, \xi_u(k), u\in W, k\in{\mathbb N}\},
\qquad
\{\zeta_u^*, \xi_u(k), u\in W, k\in{\mathbb N}\},
\qquad
\{\zeta_u^*, \xi_u^*(k), u\in W, k\in{\mathbb N}\},
\\
&
\{\zeta_u^{\star}, \xi_u^*(k), u\in W, k\in{\mathbb N}\},
\qquad
\{\zeta_u^{\star}, \xi_u^{\star}(k), u\in W, k\in{\mathbb N}\},
\qquad
\{\zeta_u^{o}, \xi_u^{o}(k), u\in W, k\in{\mathbb N}\}
\end{align*}
consists of independent random variables. We introduce random variables
\begin{align*}
&
d_i^{(2)}=\sum_{u\in W}\sum_{k=1}^{\zeta^*_u}\xi_{u}(k),
\qquad
d_i^{(3)}
=
\sum_{u\in W}\sum_{k=1}^{\zeta_{u}^*}\xi^*_{u}(k),
\\
&
d_i^{(4)}
=
\sum_{u\in W}\sum_{k=1}^{\zeta_{u}^{\star}}\xi^{\star}_{u}(k),
\qquad
d_i^{(5)}
=
\sum_{u\in W}\sum_{k=1}^{\zeta_{u}^{o}}\xi^{o}_{u}(k).
\end{align*}

When estimating 
the total variation distance between sums of random variables, say, 
\linebreak
$\sum_{k\in [m]}\eta_k=:\eta$ and
$\sum_{k\in[m]}\kappa_k=:\kappa$  we will often apply the following device. We define intermediate sums 
$\varphi_r=\sum_{k=1}^r\kappa_k+\sum_{k=r+1}^m\eta_k$
so that $\eta=\varphi_0$ and $\kappa=\varphi_m$ and 
note that $d_{TV}(\varphi_r,\varphi_{r+1})\le d_{TV}(\eta_r,\kappa_r)$.
Combining this inequality with the triangle inequality we have
\begin{align}
\label{2024-10-05+5}
d_{TV}(\eta,\kappa)
\le
\sum_{r=0}^{m-1}d_{TV}(\varphi_r,\varphi_{r+1})
\le
\sum_{r=1}^md_{TV}(\eta_r,\kappa_r).
\end{align}
Let $a,b>0$ and let $\tau,\nu$ be random variables. 
We will use the 
following inequality for the total variation distance between compound Poisson distributions
${\cal CP}(a,P_{\tau})$ and ${\cal CP}(b,P_{\nu})$, see \cite{Bloznelis_Leskela_2023} formula (25),
\begin{align}
\label{2024-10-05+6}
d_{TV}\left({\cal CP}(a,P_{\tau}),{\cal CP}(b,P_{\nu})\right)
\le 2|a-b|+bd_{TV}(P_{\tau},P_{\nu}).
\end{align}

\bigskip

Now we are ready to prove the theorem. We first  assume that the sequence $x_1,x_2,\dots$ is bounded, that is, for some $M>0$ we have $x_j\le M$ $\forall j$. The proof of the general case (where this assumption is waived) is given afterwards. The assumption implies that 
\begin{align}
\label{2024-10-04}
p_{ju}\le \frac{y_jM}{y_jM+\mu}\le \frac{y_jM}{\mu}, 
\qquad
{\text{for}}
\qquad
(j,u)\in V\times W.
\end{align}
Here the first inequality follows from the fact that the function $x\to\frac{y_jx}{y_jx+\mu}$ is increasing.

The proof consist of two parts. In the first part (steps 1 - 5 below) we show that 
$d_{TV}(d_i,d_i^{(5)})=o(1)$. In the second part (step 6 below) we show that the 
Fourier transform (characteristic function) of the distribution of 
$d_i^{(5)}$ converges to that of ${\cal CP}(y_ia_x\gamma_o,Q)$.

Let us show that $d_{TV}(d_i,d_i^{(5)})=o(1)$. We first apply triangle inequality
\begin{align}
d_{TV}\bigl(d_i,d_i^{(5)}\bigr)
\le 
d_{TV}\bigl(d_i,d_i^{(1)}\bigr)
+
\sum_{k=1}^4d_{TV}\bigl(d_i^{(k)},d_i^{(k+1)}\bigr)
\end{align}
and then show that each term on the right is of the order $o(1)$.

{\it Step 1}. Here we show that $d_{TV}(d_i,d_i^{(1)})=o(1)$.
We have, see (\ref{2024-09-30}), (\ref{2024-09-30+3}),
\begin{align}
\label{2024-10-01}
d_{TV}(d_i,d_i^{(1)})
&
\le
\PP\{d_i\not=d_i^{(1)}\}
\le \PP\{R_i\ge 1\}
\le \E R_i
\\
\nonumber
&
\le 
\sum_{\{u,v\}\subset W}p_{iu}p_{iv}\sum_{j\in V}p_{ju}p_{jv}.
\end{align}
Using (\ref{2024-10-04}) we estimate $p_{iu}p_{iv}
\le \frac{y_i^2M^2}{\mu^2}$
and  $p_{ju}p_{jv}\le \frac{y_jM}{\mu}
\frac{y_jM}{y_jM+\mu}$ and upper bound  the right side of (\ref{2024-10-01}) by
\begin{align}
\label{2024-10-04+1}
M^3\binom{m}{2}
\frac{y_i^2}{\mu^3}
S,
\qquad
{\text{where}}
\qquad
S:=\sum_{j\in V} 
y_j
\frac{y_jM}{y_jM+\mu}.
\end{align}
It remains to  show that $S=o(n)$.
To this aim we fix 
$\varepsilon\in (0,1)$ and estimate
\begin{align}
\label{2024-10-01+6}
\frac{y_jM}{y_jM+\mu}
\le
\begin{cases}
1,
\quad
\
\qquad
\qquad
 {\text{for}} \ \ y_j>\varepsilon\mu,
\\
\frac{y_jM}{\mu}
\le \varepsilon M,
\quad {\text{for}} \ \ y_j\le \varepsilon\mu.
\end{cases}
\end{align}
We have 
\begin{align}
\label{2024-10-04+4}
S
&
\le 
\varepsilon M
\sum_{j:\, y_j\le \varepsilon \mu}y_j
+
 \sum_{j:\, y_j> \varepsilon \mu}y_j
 \le
 \varepsilon M
 S_1(y)
+
 \sum_{j:\, y_j> \varepsilon \mu}y_j
 \\
 =
 \nonumber
&
 \varepsilon M n\langle y\rangle
+
n\int_{s>\varepsilon \mu}sP_{y,n}(ds).
\end{align}
Choosing $\varepsilon=\varepsilon_n\downarrow 0$ such that 
$\varepsilon_n\mu\to\infty$ as $n\to\infty$ we obtain  
$\int_{s>\varepsilon_n \mu}sP_{y,n}(ds)=o(1)$ by the uniform integrability condition (ii). Hence  $S=o(n)$.

\medskip

{\it Step 2.} 
Here we show that $d_{TV}\left(d_i^{(1)}, d_i^{(2)}\right)=o(1)$.
To this aim we write $d_i^{(1)}$ in the form
$d_i^{(1)}=\sum_{u\in W}\sum_{k=1}^{\zeta_u}\xi_{u}(k)$ and
apply (\ref{2024-10-05+5}). 
In this case $\eta_r=\sum_{k=1}^{\zeta_{r}}\xi_{r}(k)$ 
and $\kappa_r=\sum_{k=1}^{\zeta^*_{r}}\xi_{r}(k)$.
Invoking inequalities 
$d_{TV}(\eta_r,\kappa_r)\le d_{TV}(\zeta^*_{r},\zeta_{r})\le p_{ir}^2$
(the last inequality follows by LeCam's inequality  \cite{Steele_1994}) 
 we obtain
\begin{align*}
d_{TV}\left(d_i^{(1)}, d_i^{(2)}\right)
\le 
\sum_{r\in [m]}d_{TV}(\zeta_{r},\zeta^*_{r})
\le \sum_{r\in [m]}p_{ir}^2.
\end{align*}
We
note that the sum on the right $\sum_{r\in [m]}p_{ir}^2
\le M^2y_i^2\frac{m}{\mu^2}=O(n^{-1})$.

\medskip

{\it Step 3.} 
Here
we show that $d_{TV}(d_i^{(2)},d_i^{(3)})=o(1)$.
To this aim we
apply (\ref{2024-10-05+5}) with  $\eta_r=
\sum_{k=1}^{\zeta^*_{r}}\xi^*_{r}(k)$ and 
$\kappa_r=\sum_{k=1}^{\zeta^*_{r}}\xi_{r}(k)$. Invoking  
inequalities
\begin{align*}
d_{TV}(\eta_r,\kappa_r)\le p_{ir}d_{TV}(\xi^*_{r},\xi_{r})
\le
p_{ir}\sum_{j\in V\setminus\{i\}} p_{jr}^2
\end{align*}
(the first inequality follows by (\ref{2024-10-05+6}), the second one follows by LeCam's inequality) we obtain
\begin{align*}
d_{TV}(d_i^{(2)},d_i^{(3)})
\le 
\sum_{r\in[m]}d_{TV}
\left(\sum_{k=1}^{\zeta^*_{r}}\xi^*_{r}(k),
\sum_{k=1}^{\zeta^*_{r}}\xi_{r}(k)\right)
\le
\sum_{r\in[m]}
p_{ir}\sum_{j\in V\setminus\{i\}} p_{jr}^2.
\end{align*}
To show that the quantity on the right is of the order $o(1)$ we proceed as in (\ref{2024-10-01+6}),
(\ref{2024-10-04+4}) above. We have 
\begin{align*}
\sum_{r\in[m]}p_{ir}\sum_{j\in V\setminus\{i\}}p_{jr}^2
&
\le
\sum_{r\in[m]}\frac{y_iM}{\mu}
\sum_{j\in V\setminus\{i\}}\frac{y_jM}{\mu}\frac{y_jM}{y_jM+\mu}
\\
&=M^2\frac{y_i}{n}\sum_{j\in V\setminus\{i\}}y_j\frac{y_jM}{y_jM+\mu}
\le M^2\frac{y_i}{n}S=o(1).
\end{align*}

\medskip

{\it Step 4.} 
Here
we show that $d_{TV}(d_i^{(3)},d_i^{(4)})=o(1)$.
To this aim we
apply (\ref{2024-10-05+5}) with  $\eta_r=\sum_{k=1}^{\zeta_{r}^*}\xi^*_{r}(k)$ and 
$\kappa_r=\sum_{k=1}^{\zeta^{\star}_{r}}\xi^{\star}_{r}(k)$. Invoking inequalities
\begin{align*}
d_{TV}(\eta_r,\kappa_r)\le 2(p_{ir}^{\star}-p_{ir})+p_{ir}d_{TV}(\xi^{\star}_{r},\xi^*_{r})
\le
2(p_{ir}^{\star}-p_{ir})
+
2p_{ir}\sum_{j\in V\setminus\{i\}} (p^{\star}_{jr}-p_{jr})
\end{align*}
(the first inequality follows by (\ref{2024-10-05+6}), the second one follows from  (\ref{2024-10-05+6}) applied to
Poisson random variables  $\xi^{\star}_{r}$ and $\xi^*_{r}$) we obtain
\begin{align}
\label{2024-10-05+9}
d_{TV}(d_i^{(3)},d_i^{(4)})
\le 
2\sum_{r\in[m]}(p_{ir}^{\star}-p_{ir})
+
2
\sum_{r\in[m]}p_{ir}
\sum_{j\in V\setminus\{i\}}
 (p^{\star}_{jr}-p_{jr}).
\end{align}
Using inequality $0\le p^{\star}_{ir}-p_{ir}\le (p^{\star}_{ir})^2$, which 
follows from the first inequality of (\ref{2024-09-10}), we upper 
bound the first
term on the right of (\ref{2024-10-05+9})
\begin{align}
\label{2024-10-07+1}
\sum_{r\in[m]}(p_{ir}^{\star}-p_{ir})
\le 
\sum_{r\in [m]}(p_{ir}^{\star})^2
\le 
\sum_{r\in [m]}\frac{y_i^2M^2}{\mu^2}
=M^2\frac{y_i^2}{n}=O(n^{-1}).
\end{align}
To show that the second term on the right  of (\ref{2024-10-05+9}) is of 
the order $o(1)$ we proceed as in (\ref{2024-10-01+6}),
(\ref{2024-10-04+4}) above. Denote $S'=\sum_{r\in[m]}p_{ir}
\sum_{j\in V\setminus\{i\}}
 (p^{\star}_{jr}-p_{jr})$. Given $\varepsilon>0$ 
we estimate
\begin{align*}
 p^{\star}_{jr}-p_{jr}
 \le
 \begin{cases}
 p^{\star}_{jr}\le y_j\frac{M}{\mu},
 \qquad
 \qquad
 \qquad
 \quad
 \
 \
 {\text{for}}
 \quad
 y_j>\varepsilon\mu,
 \\
 (p^{\star}_{jr})^2\le \frac{y_j^2M^2}{\mu^2}\le \varepsilon y_j\frac{M^2}{\mu},
  \qquad
 {\text{for}}
 \quad
 y_j\le \varepsilon\mu.
 \end{cases}
 \end{align*}
Furthermore, we estimate $p_{ir}\le p_{ir}^{\star}$. These inequalities imply
\begin{align*}
S'
 &
 \le 
 \sum_{r\in[m]}\frac{y_iM}{\mu}
\left(
\varepsilon \frac{M^2}{\mu}
\sum_{j:\, y_j\le \varepsilon \mu} 
 y_j
 +
 \frac{M}{\mu}\sum_{j:\, y_j>\varepsilon \mu} 
 y_j
 \right)
 \\
 &
 \le
 y_iM^2\frac{nm}{\mu^2}\left(\varepsilon M\langle y\rangle
 +
 \int_{s>\varepsilon \mu}sP_{y,n}(ds)\right).
  \end{align*}
Choosing $\varepsilon=\varepsilon_n\downarrow 0$ such that 
$\varepsilon_n\mu\to\infty$ as $n\to\infty$ we obtain  
$\int_{s>\varepsilon_n \mu}sP_{y,n}(ds)=o(1)$ by the uniform integrability condition (ii). Hence $S'=o(1)$.  
 
\medskip

{\it Step 5.} 
Here
we show that $d_{TV}(d_i^{(4)},d_i^{(5)})=o(1)$.
To this aim we
apply (\ref{2024-10-05+5}) with  $\eta_r=\sum_{k=1}^{\zeta_{r}^o}\xi^o_{r}(k)$ and 
$\kappa_r=\sum_{k=1}^{\zeta^{\star}_{r}}\xi^{\star}_{r}(k)$. Invoking inequalities
\begin{align*}
d_{TV}(\eta_r,\kappa_r)
&
\le 
2\left|\frac{\gamma_o}{\gamma}-1\right|p_{ir}^{\star}
+
p_{ir}^{\star}d_{TV}(\xi^{\star}_{r},\xi^o_{r})
\\
&
\le
2\left|\frac{\gamma_o}{\gamma}-1\right|p_{ir}^{\star}
+
p_{ir}^{\star}
2
\left|\E \xi_r^{\star}-\E \xi_r^{o}\right|
\end{align*}
(the first inequality follows by (\ref{2024-10-05+6}), the second one follows from  (\ref{2024-10-05+6}) applied to
Poisson random variables  $\xi^{\star}_{r}$ and $\xi^{o}_{r}$) we obtain
\begin{align}
\label{2024-10-07}
d_{TV}(d_i^{(4)},d_i^{(5)})
\le 
2\left|\frac{\gamma_o}{\gamma}-1\right|
\sum_{r\in[m]}p_{ir}^{\star}
+
2
\sum_{r\in[m]}p_{ir}^{\star}
\left|\E \xi_r^{\star}-\E \xi_r^{o}\right|
=:R_1+2R_2.
\end{align}
The first term on the right $R_1=o(1)$ because 
$\sum_{r\in[m]}p_{ir}^{\star}=y_i\gamma\langle x\rangle\le y_i\gamma M$ is bounded and
$\gamma\to \gamma_o$. To show that $R_2=o(1)$ we write
$\E \xi_r^{o}-\E \xi_r^{\star}$ in the form
\begin{align*}
\E \xi_r^{o}-\E \xi_r^{\star}
=
x_r\delta+p_{ir}^{\star},
\qquad
\delta:=\gamma_o^{-1}a_y-\gamma^{-1}\langle y\rangle
\end{align*}
and note that $\delta=o(1)$.  We have
\begin{align*}
R_2
\le 
|\delta|\sum_{r\in[m]}x_rp_{ir}^{\star}
+
\sum_{r\in[m]}(p_{ir}^{\star})^2
=o(1).
\end{align*}
Here we used 
 $\sum_{r\in[m]}x_rp_{ir}^{\star}
\le y_iM^2\gamma=O(1)$ and 
$\sum_{r\in[m]}(p_{ir}^{\star})^2=O(n^{-1})$, see (\ref{2024-10-07+1}).

\medskip

{\it Step 6.} We write the Fourier transform of the probability distribution ${\cal CP}(y_ia_x\gamma_o,Q)$ in the form 
$f(t)=e^{y_ia_x\gamma_o(f_Q(t)-1)}$, $t\in {\mathbb R}$. Here 
$f_Q(t)=\int e^{{\bf i}ts}Q(ds)$ denotes the Fourier transform of the probability distribution $Q$; ${\bf i}=\sqrt{-1}$ denotes the imaginary unit. 
We write the characteristic function of $d_i^{(5)}$ in the form
(recal that $S_1(x)=\sum_{r=1}^mx_r$)
\begin{align*}
 \E e^{{\bf i}td_i^{(5)}}
 &
 =
\prod_{r=1}^m
\exp\left\{\left(e^{(e^{{\bf i}t}-1)\E\xi_r^o} -1\right)\E\zeta_r^o\right\}
\\
&
=
\exp\left\{
\gamma_oy_i\langle x\rangle \sum_{r=1}^m\frac{x_r}{S_1(x)}\left(e^{(e^{{\bf i}t}-1)\E\xi_r^o} -1\right)
\right\}.
\end{align*}
We denote by  ${\tilde P}_X(ds):=\frac{s}{a_x}P_X(ds)$  the size biased distribution $P_X$. We denote by 
${\tilde P}_{x,m}(ds)=\sum_{r=1}^m\frac{x_r}{S_1(x)}\delta_{x_r}$ 
 the size biased distribution $P_{x,m}$.
Condition (i) implies that ${\tilde P}_{x,m}$ converges weakly to 
${\tilde P}_X$ as $m\to+\infty$.
 Hence 
 \begin{align*}
  \sum_{r=1}^m\frac{x_r}{S_1(x)}\left(e^{(e^{{\bf i}t}-1)\E\xi_r^o} -1\right)
  =
  \int 
  \left(e^{(e^{{\bf i}t}-1)\lambda_s} -1\right)
  {\tilde P}_{x,m}(ds)
 \end{align*}
  converges to
   \begin{align*}
\int \left(e^{(e^{{\bf i}t}-1)\lambda_s} -1\right)P_{X}(ds)
=
f_Q(t)-1.
 \end{align*}
 This fact together with the convergence of the first moments $\langle x\rangle\to a_x$   yield the pointwise convergence  of the Fourier transforms
 $\E e^{{\bf i} td_i^{5)}}
\to f(t)$ as $m\to\infty$. 

We have proved Theorem \ref{T2} in the case, where 
the sequence $x_k, k\ge 1$ is bounded.

Now we waive the assumption  of boundedness of $x_k,k\ge 1$.
Let $x_r$, $ r\ge 1$ be a weight sequence satisfying condition (i).
Given $M>0$ 
define the truncated sequence 
$x_{r,M}=x_r{\mathbb I}_{\{x_r\le M\}}$, $r\ge 1$.
Let $P_{X,M}$ denote the probability distribution of $X{\mathbb I}_{\{X\le M\}}$.
Condition (i) implies that $\frac{1}{m}\sum_{r\in [m]}\delta_{x_{r,M}}$ converges weakly to $P_{X,M}$ and 
$\frac{1}{m}\sum_{r\in [m]}x_{r,M}$ converges to
 $\E (X{\mathbb I}_{\{X\le M\}})=:a_{x,M}$   as $m\to\infty$.
Let $d_{i,M}$ denote the degree of vertex $i\in V$ in the affiliation random graph $G'_M$ defined by the sequences $y_k, k\ge 1$ and $x_{k,M}, k\ge 1$. 
Let $D_M$  and $D$ be random variables with the distributions 
${\cal CP}(y_ia_{x,M}\gamma_o,Q_M)$ and ${\cal CP}(y_ia_{x}\gamma_o,Q)$.
We have
\begin{displaymath}
d_{TV}(d_i,D)
\le 
d_{TV}(d_i,d_{i,M})
+
d_{TV}(d_{i,M},D_M)
+
d_{TV}(D_M,D).
\end{displaymath}
Note that the first term on the right 
\begin{align*}
d_{TV}(d_i, d_{i,M})
&
\le 
\PP\{d_{i,M}\not= d_i\}
\le
 \PP\{\exists u\in W: i\sim u, x_u>M\} 
\\
&
\le
\E\left( \sum_{u\in W}{\mathbb I}_{\{i\sim u\}}
{\mathbb I}_{\{x_u>M\}}
\right)
=\sum_{u\in W}p_{iu}
{\mathbb I}_{\{x_u>M\}}
\\
&
\le
\sum_{u\in W}p_{iu}^{\star}
{\mathbb I}_{\{x_u>M\}}
=y_i\gamma
\int_{s>M}sP_{x,m}(ds)
\end{align*}
 converges to $0$ uniformly in $m$ as $M\to+\infty$, by the uniform integrability condition (i).
Furthermore, the second term 
$d_{TV}(d_{i,M},D_M)=o(1)$ as $n,m\to\infty$ because the convergence in distribution of integer valued random variables implies the convergence in the total variation distance.
Finally, 
 for $M\to+\infty$ we have 
$d_{TV}(D_M,D)=o(1)$, because $a_{x,M}\to a_x$ and $Q_M\to Q$.

Now we show that $d_{TV}(d_i,D)\to 0$ as $n,m\to\infty$. We fix $\varepsilon>0$ and choose large $M$ such that 
$d_{TV}(d_i,d_{i,M})<\varepsilon$ and $d_{TV}(D_M,D)<\varepsilon$.
Then, given $M$, we let $n,m\to\infty$. 
We obtain $d_{TV}(d_i,D)\le 2\varepsilon+o(1)$.
Proof of Theorem \ref{T2} is complete.

\medskip

Before the proof of Theorem \ref{T3+} below we introduce some notation and state an auxiliary result.
Let  $V_0^3$ denote the set of ordered triples $(i,j,k)$ of distinct elements $i,j,k\in V$; by  $\binom{V}{3}$ we denote the collection of subsets of $V$ of size $3$. 
For $w\in W$ we denote  
\begin{align*}
U_w=\sum_{i\in V}p_{iw}^{\star},
\qquad
T(w)=\sum_{\{i,j,k\}\in \binom{V}{3}}
 p^{\star}_{iw}p^{\star}_{jw}p^{\star}_{kw}.
 \end{align*}
 Furthermore, we denote 
 \begin{align*}
 T'(u,v)=\sum_{(i,j,k)\in V_0^3}
 p^{\star}_{iu}p^{\star}_{ju}p^{\star}_{jv} p^{\star}_{kv}
 \qquad
 {\text{and}}
 \qquad
 T'
=
 \sum_{\{u,v\}\subset W}
T'(u,v).
 \end{align*}
\begin{lem}\label{Lemma1}
Let $\mu^2=nm$.
Denote $L_{y}= \langle y^2\rangle
  \langle y\rangle^2
-
\frac{2}{n}
 \langle y^3\rangle
  \langle y\rangle
-
\frac{1}{n}
\langle y^2\rangle^2
+
\frac{2}{n^2}
 \langle y^4\rangle$.
We have
\begin{align}
\label{2024-11-04}
&
0
\le
\frac{1}{6}\frac{n^3}{\mu^3}
x_w^3
\langle y\rangle^3
-
T(w)
\le 
\frac{1}{2}
\frac{n^2}{\mu^3}
x_w^3
\langle y^2\rangle
\langle y\rangle,
\\
\label{2024-11-04+}
&
T'_3
 =
 \frac{n}{2}
\left(\langle x^2\rangle^2-\frac{\langle x^4\rangle}{m}\right)
 L_y.
\end{align}
\end{lem}
{\it Proof of Lemma \ref{Lemma1}.}
Recall that $V=[n]$. 
For numbers $a_i,b_i,c_i$, $i\in[n]$  we use 
notation $S(a)=\sum_{i\in [n]}a_i$, 
$S(ab)=\sum_{i\in [n]}a_ib_i$, $S(abc)=\sum_{i\in[n]}a_ib_ic_i$.
 
 Proof of (\ref{2024-11-04}).
We apply identity 
\begin{align*}
S^3(a)
=\sum_{i\in[n]}a_i^3+3\sum_{i\in [n]}a_i^2\sum_{j\in [n]\setminus{i}}a_j
+6\sum_{1\le i<j<k\le n}a_ia_ja_k
\end{align*}
to $a_i=p_{iw}^{\star}$ and obtain inequalities that are equivalent to (\ref{2024-11-04}).
\begin{align*}
0
\le
U_w^3
-
6T(w)
\le 
3
\sum_{i\in V}(p_{iw}^{\star})^2\sum_{j\in V}p_{jw}^{\star}.
\end{align*}
Proof of (\ref{2024-11-04+}). We apply identity
\begin{align*}
 \sum_{(i,j,k)\in V_0^3}a_ib_jc_k
 =
 S(a)S(b)S(c)
 -S(a)S(bc)
 -S(b)S(ac)
 -S(c)S(ab)
 +
 2S(abc)
 \end{align*}
 to $a_i=p_{iu}^{\star}$, 
 $b_j=p_{ju}^{\star}p_{jv}^{\star}$, $c_k=p_{kv}^{\star}$ and obtain
 $
 T'(u,v)
=
 \frac{n^3}{\mu^4}
 x_u^2x_v^2
L_y
$.

\medskip

{\it Proof of Theorem \ref{T3+}}. 
Let $G_c=(V, E_c)$ be the multigraph with colored edges 
 defined by the bipartite graph $H=(V\cup W, E)$:
vertices $i,j\in V$ are connected by an edge of color $w\in W$ in $G_c$
(denoted $i\overset {w}{\sim}j$) whenever 
$i,j$ are neighbours of $w$ in $H$. 
A subgraph of $G_c$ is projected  to subgraph of $G'$ 
 by removing edge colors and merging obtained parallel edges.

{\it Triangle count.} Here we show (\ref{2024-10-15+1}).
 Note that each triangle of $G'$ is either projection  of a monochromatic triangle or  projection of  a  triangle with all edges of different colors (we call such triangle of $G_c$  polychromatic). 
Let $\Delta_w$ be the set of monochromatic triangles of color $w$. Let
$\Delta_{u,v,z}$ be the set of polychromatic triangles with edge colors 
$u,v,z$.
 Let $\Delta_w^*$ (respectively $\Delta_{u,v,z}^*$)
 denote the set of triangles in $G'$ obtained as projections of triangles 
 from $\Delta_w$ (respectively $\Delta_{u,v,z}$).
We denote 
\begin{align*}
\Delta'=\cup_{w\in W}\Delta_w^*
\qquad
{\text{  and}}
\qquad
\Delta''=\cup_{\{u,v,z\}\subset W}\Delta_{u,v,z}^*.
\end{align*}
The union $\Delta'\cup\Delta''$ contains all triangles of  $G'$. Hence $N_{\Delta}=|\Delta'\cup\Delta''|$ and 
 (\ref{2024-10-15+1}) follows from the relations
shown below 
\begin{align}
\label{2024-10-14}
|\Delta'|
&
=\frac{n}{6\gamma}\langle x^3\rangle \langle y\rangle^3+
O_P(\sqrt{n}),
\\
\label{2024-10-14+1}
\E |\Delta''|
&
\le
\frac{1}{6}
\langle x^2\rangle^3\langle y^2\rangle^3
=
O_P(1).
\end{align}

{\it Proof of (\ref{2024-10-14}).}
Denote  $T_1=\sum_{w\in W}|\Delta_w^*|$ 
and 
$T_2=\sum_{\{u,v\}\subset W}|\Delta_u^*\cap\Delta_v^*|$.
By  inclusion-exclusion inequalities we have
\begin{align}
\label{2024-10-11}
T_1
-
T_2
\le
|\Delta'|
\le 
T_1.
\end{align}
We show  below that $\E T_2=o(1)$. Hence   
$|\Delta'|
= 
T_1
+o_P(1)$. Furthermore, we show that
\begin{align}
\label{2024-10-15+2}
\E T_1
=
\frac{n}{6\gamma}
\langle x^3\rangle
\langle y\rangle^3
+O(1)
\end{align}
and $\Var T_1=O(n)$. 
These two relations imply $T_1=\E T_1+O_P(\sqrt{n})$, by 
Chebyshev's inequality.  We have arrived to (\ref{2024-10-14}).

Let us show that  $\E T_2=O(m^{-1})$. We have 
$|\Delta_u^*\cap\Delta_v^*|
=
\sum_{\{i,j,k\}\subset V}
{\mathbb I}_{iu}
{\mathbb I}_{ju}
{\mathbb I}_{ku}
{\mathbb I}_{iv}
{\mathbb I}_{jv}
{\mathbb I}_{kv}$
and
\begin{align}
\nonumber
\E T_2
&
=
\E\left(\sum_{\{u,v\}\subset W}
\sum_{\{i,j,k\}\subset V}
{\mathbb I}_{iu}
{\mathbb I}_{ju}
{\mathbb I}_{ku}
{\mathbb I}_{iv}
{\mathbb I}_{jv}
{\mathbb I}_{kv}
\right)
\\
&
\le
\sum_{\{u,v\}\subset W}
\sum_{\{i,j,k\}\subset V}
p_{iu}^{\star}
p_{ju}^{\star}
p_{ku}^{\star}
p_{iv}^{\star}
p_{jv}^{\star}
p_{kv}^{\star}
\\
\label{2024-10-15+23}
&
\le
\frac{1}{12m}\langle x^3\rangle^2\langle y^2\rangle^3.
\end{align}

Let us prove (\ref{2024-10-15+2}). Note that $T_1=\sum_{w\in W}|\Delta_w|$ since 
$|\Delta_w^*|=|\Delta_w|$ $\forall w\in W$.  We evaluate \begin{align*}
\E|\Delta_w|
=
\sum_{\{i,j,k\}\subset V}
{\mathbb I}_{iw}
{\mathbb I}_{jw}
{\mathbb I}_{kw}
=\sum_{\{i,j,k\}\subset V}p_{iw}p_{jw}p_{kw}.
\end{align*}
 Using  (\ref {2024-09-10}) we approximate 
 $\E |\Delta_w|$ by $T(w):=\sum_{\{i,j,k\}\subset V}
 p^{\star}_{iw}p^{\star}_{jw}p^{\star}_{kw}$, 
\begin{align}
\label{2024-10-14+3}
&
T(w)-R_1(w)
\le
\E|\Delta_w|
\le T(w),
\\
\nonumber
&
R_1(w)
:=
\sum_{\{i,j,k\}\subset V}
p^{\star}_{iw}p^{\star}_{jw}p^{\star}_{kw}
\left(p^{\star}_{iw}+p^{\star}_{jw}+p^{\star}_{kw}\right).
\end{align}
A straightforward calculation shows
\begin{align*}
R_1(w)
=\frac{1}{6}
\sum_{(i,j,k)\in V_0^3}
p^{\star}_{iw}p^{\star}_{jw}p^{\star}_{kw}
\left(p^{\star}_{iw}+p^{\star}_{jw}+p^{\star}_{kw}\right) 
\le 
\frac{1}{2}\frac{n^3}{\mu^4}
x_w^4
\langle y^2\rangle
\langle y\rangle^2.
\end{align*}

Combining these inequalities with (\ref{2024-11-04}) and (\ref{2024-10-14+3}) 
we obtain
\begin{align}
0
\le 
\frac{1}{6}\frac{n^3}{\mu^3}
x_w^3
\langle y\rangle^3
-
\E|\Delta_w|
\le 
\frac{1}{2}\frac{n^3}{\mu^4}
x_w^4
\langle y^2\rangle
\langle y\rangle^2
+
 \frac{1}{2}
 \frac{n^2}{\mu^3}
 x_w^3
 \langle y^2\rangle
 \langle y\rangle.
 \end{align}
Summing over $w\in W$  we obtain
for $\mu=\sqrt{nm}$
\begin{align}
0
\le 
\frac{n}{6\gamma}\langle x^3\rangle\langle y\rangle^3
-
\E T_1
\le 
\frac{1}{2\gamma^2}
\langle x^4\rangle\langle y^2\rangle\langle y\rangle^2
+
\frac{1}{2\gamma}
\langle x^3\rangle\langle y^2\rangle\langle y\rangle.
\end{align}
We have arrived to (\ref{2024-10-15+2}).

It remains to show that $\Var \, T_1=O(n)$.
By the independence of $\Delta_w$, $w\in W$, we have
$\Var \, T_1=\sum_{w\in W}\Var |\Delta_w|$. 
We show below  that for each $w\in W$
\begin{align}
\label{2024-10-18}
\Var |\Delta_w|
\le
 \sum_{r=3}^5\frac{x_w^r}{\gamma^r}\langle y\rangle^r.
\end{align}
This bound implies
\begin{align*}
\Var\, T_1\le \sum_{r=3}^5\frac{m}{\gamma^r}\langle x^r\rangle\langle y\rangle^r
=
n\sum_{r=3}^5\gamma^{2-r}\langle x^r\rangle\langle y\rangle^r
=
O(n).
\end{align*}

Proof of (\ref{2024-10-18}). We fix $w\in W$.
For $A \subset V$ we denote
${\mathbb I}_{A}=\prod_{i\in A}{\mathbb I}_{iw}$ and 
${\bar {\mathbb I}}_{A}={\mathbb I}_{A}-\E {\mathbb I}_{A}$.
Note that $|\Delta_w|=\sum_{A\in\binom{V}{3}}{\mathbb I}_{A}$.
Since
$\E ({\bar {\mathbb I}}_{A}{\bar {\mathbb I}}_{A'})=0$ for $A\cap A'=\emptyset$, we have
\begin{align}
\label{2024-10-15}
&\Var |\Delta_w|
=
\E\left(\sum_{A\in \binom{V}{3}}{\bar {\mathbb I}}_A\right)^2
=
S_1+2S_2+2S_3,
\\
\nonumber
&
S_1:=
\sum_{A\in \binom{V}{3}}\E{\bar {\mathbb I}}_A^2,
\qquad
S_2:=\sum_{\{A, A'\}\subset \binom{V}{3}: |A\cap A'|=2}
\E({\bar {\mathbb I}}_A{\bar {\mathbb I}}_{A'}),
\\
\nonumber
&
S_3:=\sum_{\{A, A'\}\subset \binom{V}{3}: |A\cap A'|=1}
\E({\bar {\mathbb I}}_A{\bar {\mathbb I}}_{A'}).
\end{align}
 We upper bound the sums $S_1,S_1,S_3$  using simple inequality 
 $\E({\bar {\mathbb I}}_A{\bar {\mathbb I}}_{A'})
 \le \E {\mathbb I}_{A\cup A'}$.
 A straightforward calculation shows that
 \begin{align*}
 S_1
 &
 \le
 \sum_{A\in\binom{V}{3}}\E{\mathbb I}_A
=
 \sum_{\{i,j,k\}\subset V}p_{iw}p_{jw}p_{kw}
 \le 
 \sum_{\{i,j,k\}\subset V}p_{iw}^{\star}p_{jw}^{\star}p_{kw}^{\star}
\le
\frac{U_w^3}{6}
=
 \frac{x_w^3}{6\gamma^3}\langle y\rangle^3,
 \\
 S_2
 &
 \le
\sum_{\{i,j\}\subset V}
\
\sum_{\{k,l\}\subset V\setminus\{i,j\}}\E {\mathbb I}_{\{i,j,k\}\cup\{i,j,l\}}=
\sum_{\{i,j\}\subset V}
\
\sum_{\{k,l\}\subset V\setminus\{i,j\}}p_{iw}p_{jw}p_{kw}p_{lw}
\\
&\le 
\sum_{\{i,j\}\subset V}
\
\sum_{\{k,l\}\subset V\setminus\{i,j\}}
p_{iw}^{\star}p_{jw}^{\star}p_{kw}^{\star}p_{lw}^{\star}
\le 
\frac{U_w^4}{4}
=
\frac{x_w^4}{4\gamma^4}\langle y\rangle^4,
\\
 2S_3
 &
 \le
\sum_{i\in V}
\
\sum_{\{j,k\}\subset V\setminus\{i\}}
\
\sum_{\{s,t\}\subset V\setminus\{i,j,k\}}\E {\mathbb I}_{\{i,j,k\}\cup\{i,s,t\}}
\\
&
 =
\sum_{i\in V}
\
\sum_{\{j,k\}\subset V\setminus\{i\}}
\
\sum_{\{s,t\}\subset V\setminus\{i,j,k\}}p_{iw}p_{jw}p_{kw}p_{sw}p_{tw}
\\
&
\le
\frac{1}{4}
\left(\sum_{i\in V}p_{iw}\right)^5
\le
\frac{U_w^5}{4}
=
\frac{x_w^5}{4\gamma^5}\langle y\rangle^5.
\end{align*}
Invoking these bounds in (\ref{2024-10-15}) and then summing
(\ref{2024-10-15}) over $w\in W$ we obtain (\ref{2024-10-18}).

{\it Proof of (\ref{2024-10-14+1}).}
Given $\{i,j,k\}\subset V$ and $\{u,v,w\}\subset W$
 the expected number of polychromatic triangles in $G_c$
  on vertices $i,j,k$ and with edge colors $u,v,w$ is the  sum  
\begin{align*}
\sum_{\pi}
p_{i\pi_u}p_{j\pi_u}p_{i\pi_v}p_{k\pi_v}p_{j\pi_w}p_{k\pi_w}
\end{align*}
that runs over permutations $\pi=(\pi_u,\pi_v,\pi_w)$ of $(u,v,w)$.
Using $p_{ab}\le p_{ab}^{\star}$ we upper bound the sum by
\begin{align*}
\sum_{\pi}p_{i\pi_u}^{\star}p_{j\pi_u}^{\star}p_{i\pi_v}^{\star}p_{k\pi_v}^{\star}p_{j\pi_w}^{\star}p_{k\pi_w}^{\star}
=\frac{6}{\mu^6} (x_ux_vx_wy_iy_jy_k)^2.
\end{align*}
Consequently, we have
\begin{align}
\E |\Delta''|
\le
\sum_{\{i,j,k\}\subset V}
\sum_{\{u,v,w\}\subset W}
\frac{6}{\mu^6} (x_ux_vx_wy_iy_jy_k)^2
\le
\frac{1}{6}
\langle x^2\rangle^3
\langle y^2\rangle^3.
\end{align}

\medskip

{\it Count of $2$-paths.} Here we show (\ref{2024-10-26+30}).
For $u,v\in W$ we denote by 
$\Lambda_{u}$ (respectively $\Lambda_{uv}$) the set of $2$-paths of 
$G_c$  with both edges colored $u$ (respectively with edges receiving 
different colors $u$ and $v$).
Every $2$-path of $G_c$ is projected to a path in $G'$ by removing edge 
colors. Let  $\Lambda_{u}^*$ (respectively $\Lambda_{uv}^*$) 
denote the sets of $2$-paths in $G'$ obtained as projections of  
$2$-paths  of $\Lambda_{u}$ (respectively $\Lambda_{uv}$).
Denote $\Lambda'=\cup_{u\in W}\Lambda_{u}^*$ and 
$\Lambda''=\cup_{\{u,v\}\subset  W}\Lambda_{uv}^*$. The union  $\Lambda'\cup\Lambda''$ contains all $2$-paths of $G'$.
Hence $N_{\Lambda}=|\Lambda'\cup\Lambda''|$. 
To show (\ref{2024-10-26+30}) we use inclusion-exclusion identity
\begin{align*}
N_{\Lambda}=|\Lambda'\cup\Lambda''|=|\Lambda'|+|\Lambda''|-
|\Lambda'\cap\Lambda''|
\end{align*}
and evaluate each term on the right
\begin{align}
\label{2024-11-01}
|\Lambda'|
&
=
3|\Delta'|
=
\frac{n}{2\gamma}\langle x^3\rangle \langle y\rangle^3+
O_P(\sqrt{n}),
\\
\label{2024-10-15+11}
|\Lambda''|
&
=
\frac{n}{2}
\langle x^2\rangle^2
\langle y\rangle^2
\langle y^2\rangle+O_P(\sqrt{n}).
\\
\label{2024-10-15+10}
\E |\Lambda'\cap\Lambda''|
&\le
\frac{2}{\gamma}
\langle x^3\rangle
\langle x^2\rangle
\langle y^2\rangle^2
\langle y\rangle
=O(1).
\end{align}
The first identity of (\ref{2024-11-01}) is obvious, the second relation
is shown in  (\ref{2024-10-14}). 
It remains to prove (\ref{2024-10-15+11}), (\ref{2024-10-15+10}).

{\it Proof of (\ref{2024-10-15+10}). } 
Given a path $i\sim j\sim k\in \Lambda'\cap \Lambda''$ there exists $\{u,v\}\in V$ such that the monochromatic $2$-path $i\overset {u}{\sim}j\overset {u}{\sim} k$ of color $u$ is present in $G_c$ and at least one of the edges 
$i\overset {v}{\sim}j$, $j\overset {v}{\sim} k$  of color $v$ is present in $G_c$. In the inequality below 
${\mathbb I}_{iu}{\mathbb I}_{ju}{\mathbb I}_{ku}$ represents the indicator of the monochromatic triangle and 
${\mathbb I}_{iv}{\mathbb I}_{jv}$ represents the indicator of  additional edge of color $v$ connecting  $i$ and $j$ in $G_c$.
 We have
  \begin{align*}
  \E |\Lambda'\cap\Lambda''|
 &
  \le
  \E\left(\sum_{u\in W}
  \sum_{v\in W\setminus \{u\}}
  \sum_{\{i,j,k\}\subset V}
  2
  {\mathbb I}_{iu}{\mathbb I}_{ju}{\mathbb I}_{ku}
  \left(
  {\mathbb I}_{iv}{\mathbb I}_{jv}
  +
  {\mathbb I}_{iv}{\mathbb I}_{kv}+{\mathbb I}_{jv}{\mathbb I}_{kv}
  \right)
  \right)
  \\
  &
  =
 2
  \sum_{u\in W}
  \sum_{v\in W\setminus \{u\}}
  \sum_{\{i,j,k\}\subset V}
  p_{iu}p_{ju}p_{ku}
  \left(
  p_{iv}p_{jv}
  +
  p_{iv}p_{kv}
  +
  p_{jv}p_{kv}
  \right)
  \\
  &
  \le
 2
  \sum_{u\in W}
  \sum_{v\in W\setminus\{u\}}
  \sum_{(i,j,k)\in V^3}
  p_{iu}p_{ju}p_{ku}p_{iv}p_{jv}
  \\
  &
  \le
 2
  \sum_{u\in W}
  \sum_{v\in W\setminus\{u\}}
  \sum_{(i,j,k)\in V^3}
  p_{iu}^{\star}p_{ju}^{\star}p_{ku}^{\star}p_{iv}^{\star}p_{jv}^{\star}
\le
\frac{2}{\gamma}
\langle x^3\rangle
\langle x^2\rangle
\langle y^2\rangle^2
\langle y\rangle.
 \end{align*} 
  
  {\it Proof of (\ref{2024-10-15+11}). } In the proof we use relations
  \begin{align}
  \label{2024-10-22}
 &
  |\Lambda_{uv}^*|\le |\Lambda_{uv}|
  =
   \sum_{(i,j,k)\in V_0^3}
 {\mathbb I}_{iu}
 {\mathbb I}_{ju} {\mathbb I}_{jv} {\mathbb I}_{kv},
 \\
 \label{2024-10-22+1}
 &
 |\Lambda_{uv}|-|\Lambda_{uv}^*|=|\Lambda_u^*\cap  \Lambda_v^*|=
 3|\Delta_u^*\cap \Delta_v^*|.
 \end{align}
We only comment on (\ref{2024-10-22+1}) since   (\ref{2024-10-22}) is obvious. The difference $|\Lambda_{uv}|- |\Lambda_{uv}^*|$ counts $2$-paths $i\sim j\sim k \in\Lambda_{uv}^*$ such that both colored paths 
$i\overset {u}{\sim}j \overset {v}{\sim}k$ and $i\overset {v}{\sim}j\overset {u}{\sim}k$ are present in $\Lambda_{uv}$. This means that 
  $i,j,k\in N_u\cap N_v$.
 In particular,
 we have 
 $|\Lambda_{uv}|- |\Lambda_{uv}^*|
 =|\Lambda_u^*\cap \Lambda_v^*|$.  Hence the first identity of  
 (\ref{2024-10-22+1}).
 The second identity of (\ref{2024-10-22+1}) follows from the fact that in $G'$  a 
 triangle  gives rise to three $2$-paths.

 Denote $T_3=\sum_{\{u,v\}\subset W}|\Lambda_{uv}|$.
 We establish (\ref{2024-10-15+11}) in three steps: we show that
 \begin{align}
 \label{2024-11-01+1}
 &
 |\Lambda''|
 =
 \sum_{\{u,v\}\subset W}|\Lambda_{uv}^*|+O_P(1),
 \\
 &
 \sum_{\{u,v\}\subset W}|\Lambda_{uv}^*|=T_3+O_P(1),
 \\
  \label{2024-10-22+2}
 & 
 T_3
  =\frac{n}{2}
\langle x^2\rangle^2
\langle y\rangle^2
\langle y^2\rangle+O_P(\sqrt{n}).
 \end{align}
 {\it Step 1.} Here  we prove (\ref{2024-11-01+1}).
 To this aim we apply inclusion-exclusion inequality
  \begin{align} \label{2024-10-15+20}
  0\le \sum_{\{u,v\}\subset W}|\Lambda_{uv}^*|-|\Lambda''|\le R_2,
  \qquad
  R_2
  :=
  \sum_{\substack{\{u,v\}, \{s,t\}\subset W: \\ \{u,v\}\not=\{s,t\}}}
  |\Lambda_{uv}^*\cap \Lambda_{st}^*|
  \end{align}
  and show that
  \begin{align*}
  \E R_2
  \le 
   \frac{1}{n}
 \langle x^2\rangle^4
\langle y^4\rangle\langle y^2\rangle^2
+
 \langle x^2\rangle^3
\langle y\rangle
\langle y^2\rangle
\langle y^3\rangle
 +
\frac{1}{\mu} \langle x^2\rangle^2
\langle x^3\rangle
\langle y^2\rangle^2
\langle y^3\rangle.
\end{align*}
To estimate $\E R_2$
we split
 \begin{align*}
 R_2=R_{2.0}+R_{2.1},
 \qquad
 {\text{where}}
 \qquad
 R_{2.k}=
  \sum_{\substack{\{u,v\}, \{s,t\}\subset W: \\ |\{u,v\}\cap\{s,t\}|=k}}
  |\Lambda_{uv}^*\cap \Lambda_{st}^*|
 \end{align*}
 and upper bound $\E R_{2.0}$, $\E R_{2.1}$.
 We start with  $\E R_{2.0}$.
For $\{u,v\}\cap \{s,t\}=\emptyset$ we have
 \begin{align*}
 \E
  |\Lambda_{uv}^*\cap \Lambda_{st}^*|
 &
  \le
 \E \sum_{(i,j,k)\in V_0^3}
 {\mathbb I}_{iu}
 {\mathbb I}_{ju} {\mathbb I}_{jv} {\mathbb I}_{kv} 
 \left(
 {\mathbb I}_{is}
 {\mathbb I}_{js} {\mathbb I}_{jt} {\mathbb I}_{kt} 
 +
 {\mathbb I}_{it}
 {\mathbb I}_{jt} {\mathbb I}_{js} {\mathbb I}_{ks}
 \right)
\\
 &
 =
 \E \sum_{(i,j,k)\in V_0^3}
 p_{iu}
 p_{ju} 
p_{jv} 
p_{kv} 
 \left(
p_{is}
 p_{js} 
 p_{jt} 
 p_{kt} 
 +
p_{it}
p_{jt} 
p_{js} 
p_{ks}
 \right)
 \\
 &
 \le 
 2\frac{n^3}{\mu^8}
 x_u^2x_v^2x_s^2x_t^2
 \langle y^4\rangle
 \langle y^2\rangle^2
.
 \end{align*}
Here the sum $\sum_{(i,j,k)\in V_0^3}
 {\mathbb I}_{iu}
 {\mathbb I}_{ju} {\mathbb I}_{jv} {\mathbb I}_{kv} $ runs over colored $2$-paths $i\overset {u}{\sim}j \overset {v}{\sim}k\in \Lambda_{uv}$ and 
 ${\mathbb I}_{is}
 {\mathbb I}_{js} {\mathbb I}_{jt} {\mathbb I}_{kt} 
 +
 {\mathbb I}_{it}
 {\mathbb I}_{jt} {\mathbb I}_{js} {\mathbb I}_{ks}
 $
 accounts for
 the
 matching $2$-paths from $\Lambda_{st}$.
 In the very last inequality we used $p_{pq}\le p_{pq}^{\star}$, see (\ref{2024-09-10}).
 We obtain
 \begin{align*}
\E R_{2.0}
 \le
 \sum_{\substack{\{u,v\}, \{s,t\}\subset W: \\ |\{u,v\}\cap\{s,t\}|=0}}
2\frac{n^3}{\mu^8}
x_u^2x_v^2x_s^2x_t^2
 \langle y^4\rangle\langle y^2\rangle^2
 \le 
 \frac{1}{n}
 \langle x^2\rangle^4
\langle y^4\rangle\langle y^2\rangle^2.
\end{align*}
We similarly upper bound $\E R_{2.1}$.  For  $u=s$ we have
 \begin{align*}
 \E
  |\Lambda_{uv}^*\cap \Lambda_{ut}^*|
 &
  \le
 \E \sum_{(i,j,k)\in V_0^3}
 {\mathbb I}_{iu}
 {\mathbb I}_{ju} {\mathbb I}_{jv} {\mathbb I}_{kv} 
 \left(
 {\mathbb I}_{iu}
 {\mathbb I}_{ju} {\mathbb I}_{jt} {\mathbb I}_{kt} 
 +
 {\mathbb I}_{it}
 {\mathbb I}_{jt} {\mathbb I}_{ju} {\mathbb I}_{ku}
 \right)
\\
 &
 =
 \E \sum_{(i,j,k)\in V_0^3}
 p_{iu}
 p_{ju} 
p_{jv} 
p_{kv} 
 \left( 
 p_{jt} 
 p_{kt} 
 +
p_{it}
p_{jt} 
p_{ku}
 \right)
 \\
&
\le
\frac{n^3}{\mu^6}x_u^2x_v^2x_t^2
\langle y\rangle
\langle y^2\rangle
\langle y^3\rangle
+
 \frac{n^3}{\mu^7}x_u^3x_v^2x_t^2
\langle y^2\rangle^2
\langle y^3\rangle.
 \end{align*}
Hence 
\begin{align*}
\E R_{2.1}
 &
 \le
 \sum_{\{u,v\}\subset W}
 \sum_{t\in W\setminus\{u,v\}}
 \E\left(|\Lambda_{uv}^*\cap\Lambda_{ut}^*|
 +
 |\Lambda_{uv}^*\cap\Lambda_{vt}^*|\right)
 \\
 &
 \le
 \langle x^2\rangle^3
\langle y\rangle
\langle y^2\rangle
\langle y^3\rangle
 +
\frac{1}{\mu} \langle x^2\rangle^2
\langle x^3\rangle
\langle y^2\rangle^2
\langle y^3\rangle.
\end{align*}

  {\it Step 2.}  Here we show (\ref{2024-10-22+2}).
 From (\ref{2024-10-22+1}) and (\ref{2024-10-15+23}) we obtain
 \begin{align}
 \label{2024-10-15+21}
T_3-\sum_{\{u,v\}\subset W}|\Lambda_{uv}^*|
=
3\sum_{\{u,v\}\subset W}|\Delta_u^*\cap \Delta_v^*|
=3T_2
 =o_P(1).
\end{align}

  {\it Step 3.} 
Here we prove (\ref{2024-10-22+2}).
To this aim
we show  that 
\begin{align}
\label{2024-10-22+5}
\E T_3
=
\frac{n}{2}
\langle x^2\rangle^2
\langle y\rangle^2
\langle y^2\rangle+O(1)
\qquad
{\text{and}}
\qquad
\Var\, T_3=O(n).
\end{align} 
Note that (\ref{2024-10-22+5}) combined with  Chebyshev's inequality implies
$T_3=\E T_3+O_P(\sqrt{\Var T_3})$.
 
Let us evaluate $\E T_3$. We have
$\E T_3=\sum_{\{u,v\}\subset W}
\sum_{(i,j,k)\in V_0^3}
 p_{iu}p_{ju}p_{jv}p_{kv}$,  by the second relation of (\ref{2024-10-22}). 
 We approximate $\E T_3$ by $T'
=
 \sum_{\{u,v\}\subset W}
\sum_{(i,j,k)\in V_0^3}
 p^{\star}_{iu}p^{\star}_{ju}p^{\star}_{jv} p^{\star}_{kv}$.
 Inequalities
 (\ref {2024-09-10}) imply
$0\le T'-\E T_3\le R'$, 
where
\begin{align*}
 R'
 &
 = \sum_{\{u,v\}\subset W}
\sum_{(i,j,k)\in V_0^3}
 p^{\star}_{iu}p^{\star}_{ju}p^{\star}_{jv} p^{\star}_{kv}
 (p^{\star}_{iu}+p^{\star}_{ju}+p^{\star}_{jv} +p^{\star}_{kv})
 \\
 &
 \le
 \frac{1}{\gamma}
\langle x^3\rangle
\langle x^2\rangle
\left(
\langle y^2\rangle^2
\langle y\rangle
+
\langle y^3\rangle
\langle y\rangle^2
\right)
 .
 \end{align*}
Combining this bound with  (\ref{2024-11-04+}) we obtain  the first relation of (\ref{2024-10-22+5}).

It remains to prove the second relation of (\ref{2024-10-22+5}). 
 Note that $|\Lambda_{uv}|$ and $|\Lambda_{s,t}|$ are 
 independent for $\{u,v\}\cap \{s,t\}=\emptyset$.
 Hence
\begin{align*}
\Var\, T_3
=
S_1+S_2,
\qquad 
S_k
=
  \sum_{\substack{\{u,v\}, \{s,t\}\subset W: \\ |\{u,v\}\cap\{s,t\}|=k}}
  \Cov (|\Lambda_{uv}|,|\Lambda_{st}|),
  \qquad
  k=1,2.
  \end{align*}
 In the rest part of the proof we show that $S_1=O(n)$ and $S_2=O(n)$.
 
 \medskip
Let us evaluate $S_2=\sum_{\{u,v\}\subset W}\Var |\Lambda_{uv}|$.
Given $\{u,v\}\subset W$  and $(i,j,k)\in V_0^3$ we denote 
\begin{align}
\label{2024-10-25+30}
&
{\mathbb I}_{(i,j,k)}
=
{\mathbb I}_{iu}
 {\mathbb I}_{ju}
{\mathbb I}_{jv}
{\mathbb I}_{kv},
\qquad
\qquad
\quad
{\bar {\mathbb I}}_{(i,j,k)}
=
{\mathbb I}_{(i,j,k)}
-
\E {\mathbb I}_{(i,j,k)},
\\
\nonumber
&
{\mathbb I}_{\{i,j,k\}}
=
\sum_{\pi}
{\mathbb I}_{(\pi_i,\pi_j,\pi_k)},
\qquad
\qquad
{\bar {\mathbb I}}_{\{i,j,k\}}
=
{\mathbb I}_{\{i,j,k\}}-\E {\mathbb I}_{\{i,j,k\}}.
\end{align}
Here the sum run over permutations $\pi=(\pi_i,\pi_j,\pi_k)$ of $i,j,k$.
We write $|\Lambda_{uv}|-\E |\Lambda_{uv}|$ in the form
$\sum_{A\subset V, |A|=3}{\bar {\mathbb I}}_{A}$, 
where 
${\bar {\mathbb I}}_{A}={\bar {\mathbb I}}_{\{i,j,k\}}$
for $A=\{i,j,k\}$, and note that
 ${\bar {\mathbb I}}_{A}$ and ${\bar {\mathbb I}}_{A'}$ are 
 uncorrelated unless $|A\cap A'|\ge 1$. Hence
\begin{align*}
&
 \Var |\Lambda_{uv}|
 =
\E \left(
 \sum_{A\subset V, |A|=3}{\bar {\mathbb I}}_{A}\right)^2
 =
 S_{2.1}(u,v)+S_{2.2}(u,v)+S_{2.3}(u,v),
\\
&
S_{2.k}(u,v)
=
  \sum_{\substack{A\subset V:\\ |A|=3}}
  \
  \
   \sum_{\substack{A'\subset V:\, |A'|=3,\\ |A\cap A'|=k}}
  \E({\bar {\mathbb I}}_{A}{\bar {\mathbb I}}_{A'}),
  \qquad
  k=1,2,3.
\end{align*}
Now we upper bound the terms $ S_{2.r}(u,v)$, $r=1,2,3$.
We use the simple inequality for  the covariance of non-negative random variables  
$X,Y$,
\begin{align}
\label{covariance}
\E (XY)-\E X\E Y\le \E (XY).
\end{align}
 We apply (\ref{covariance})  to sums of Bernoulli random variables.
We have for distinct $i,j,k,q,r\in V$
\begin{align}
\label{2024-10-22+10}
&
\E{\bar {\mathbb I}}^2_{(i,j,k)}
\le 
\E{\mathbb I}_{(i,j,k)}
=
p_{iu}
p_{ju}
p_{jv}
p_{kv} 
\le
p_{iu}^{\star}
p_{ju}^{\star}
p_{jv}^{\star}
p_{kv}^{\star}=
\frac{1}{\mu^4}x_u^2x_v^2y_iy_j^2y_k,
\\
\nonumber
&\E\left(
{\bar {\mathbb I}}_{(i,j,k)}
\left(
{\bar {\mathbb I}}_{\{i,j,q\}}
+
{\bar {\mathbb I}}_{\{i,q,k\}}
+{\bar {\mathbb I}}_{\{q,j,k\}}
\right)
\right)
\le
\E\left(
{\mathbb I}_{(i,j,k)}
\left(
{\mathbb I}_{\{i,j,q\}}
+
{\mathbb I}_{\{i,q,k\}}
+
{\mathbb I}_{\{q,j,k\}}
\right)
\right)
\\
\label{2024-10-25}
&
\
\
\le
\E\left(
{\mathbb I}_{(i,j,k)}
9 \left(
{\mathbb I}_{qu}
+
{\mathbb I}_{qv}
\right)
\right)
=
9
p_{iu}
p_{ju}
p_{jv}
p_{kv} 
\left(
p_{qu}+p_{qv}
\right)
\\
\label{2024-10-22+11}
&
\
\
\le
9
p_{iu}^{\star}
p_{ju}^{\star}
p_{jv}^{\star}
p_{kv}^{\star} 
\left(
p_{qu}^{\star}+p_{qv}^{\star}
\right)
=
\frac{9}{\mu^5}
x_u^2x_v^2(x_u+x_v)y_iy_j^2y_ky_q,
\\
\nonumber
&\E\left(
{\bar {\mathbb I}}_{(i,j,k)}
\left(
{\bar {\mathbb I}}_{\{i,r,q\}}
+
{\bar {\mathbb I}}_{\{r,j,q\}}
+
{\bar {\mathbb I}}_{\{r,q,k\}}
\right)
\right)
\le
\E\left(
{\mathbb I}_{(i,j,k)}
\left(
{\mathbb I}_{\{i,r,q\}}
+
{\mathbb I}_{\{r,j,q\}}
+
{\mathbb I}_{\{r,q,k\}}
\right)
\right)
\\
\label{2024-10-25+1}
&
\
\
\le
\E\left(
{\mathbb I}_{(i,j,k)}
9 \left(
{\mathbb I}_{qu}{\mathbb I}_{rv}+{\mathbb I}_{qv}{\mathbb I}_{ru}
\right)
\right)
=
9
p_{iu}
p_{ju}
p_{jv}
p_{kv} 
\left(
p_{qu}p_{rv}+p_{ru}p_{qv}
\right)
\\
\label{2024-10-22+12}
&
\
\
\le
9
p_{iu}^{\star}
p_{ju}^{\star}
p_{jv}^{\star}
p_{kv}^{\star} 
\left(
p_{qu}^{\star}p_{rv}^{\star}+p_{ru}^{\star}p_{qv}^{\star}
\right)
=
\frac{18}{\mu^6}
x_u^3x_v^3y_iy_j^2y_ky_qy_r.
\end{align}
In the first inequality of (\ref{2024-10-25})  (respectively (\ref{2024-10-25+1}))   we used  inequality
$
{\mathbb I}_{\{l,r,q\}}
\le 3({\mathbb I}_{qu}+{\mathbb I}_{qv})
$
(respectively
$
{\mathbb I}_{\{l,r,q\}}
\le 3({\mathbb I}_{qu}{\mathbb I}_{rv}+{\mathbb I}_{qv}{\mathbb I}_{ru})
$)
 which hold for any distinct $l,q,r\in V$.

Estimation of  $S_{2.3}(u,v)$. Combining
 the inequality $(\sum_{i=1}^6a_i)^2\le6\sum_{i=1}^6a_i^2$, which follows by Hölder's inequality, with (\ref{2024-10-22+10})  we obtain
\begin{align*}
\E {\bar {\mathbb I}}_{\{i,j,k\}}^2
\le
6
\sum_{\pi}
\E {\bar {\mathbb I}}^2_{(\pi_i,\pi_j,\pi_k)}
\le
\frac{6}{\mu^4}x_u^2x_v^2y_iy_jy_k(2y_i+2y_j+2y_k).
\end{align*}
Summing over $\{i,j,k\}\subset V$ we have
\begin{align}
\label{2024-10-25+10}
S_{2.3}(u,v)
=
\sum_{\{i,j,k\}\subset V}
\E {\bar {\mathbb I}}_{\{i,j,k\}}^2
\le
6\frac{n}{m^2}
x_u^2x_v^2
\langle y^2\rangle
\langle y\rangle^2.
\end{align}

Estimation of  $S_{2.2}(u,v)$. Here we use inequality (\ref{2024-10-22+11}).
We have
\begin{align}
\label{2024-10-25+11}
S_{2.2}(u,v)
&
=
\sum_{(i,j,k)\in V_0^3}
\
\
\sum_{q\in V\setminus\{i,j,k\}}
\E\left(
{\bar {\mathbb I}}_{(i,j,k)}
\left(
{\bar {\mathbb I}}_{\{i,j,q\}}
+
{\bar {\mathbb I}}_{\{i,q,k\}}
+
{\bar {\mathbb I}}_{\{q,j,k\}}
\right)
\right)
\\
\nonumber
&
\le
\frac{9}{\mu^5}
\sum_{(i,j,k)\in V_0^3}
\
\
\sum_{q\in V\setminus\{i,j,k\}}
x_u^2x_v^2(x_u+x_v)y_iy_j^2y_ky_q
\\
\nonumber
&
\le 
\frac{9}{\gamma}\frac{n}{m^2} x_u^2x_v^2(x_u+x_v)
\langle y^2\rangle
\langle y\rangle^3.
\end{align}
Estimation of  $S_{2.1}(u,v)$. Here we use inequality (\ref{2024-10-22+12}).
We have
\begin{align}
\label{2024-10-25+12}
S_{2.1}(u,v)
&
=
\sum_{(i,j,k)\in V_0^3}
\
\
\sum_{\{q,r\}\subset
V\setminus\{i,j,k\}}
\E\left(
{\bar {\mathbb I}}_{(i,j,k)}
\left(
{\bar {\mathbb I}}_{\{i,r,q\}}
+
{\bar {\mathbb I}}_{\{r,j,q\}}
+
{\bar {\mathbb I}}_{\{r,q,k\}}
\right)
\right)
\\
\nonumber
&
\le
\frac{18}{\mu^6}
\sum_{(i,j,k)\in V_0^3}
\
\
\sum_{\{q,r\}\subset
V\setminus\{i,j,k\}}
x_u^3x_v^3y_iy_j^2y_ky_qy_r
\\
\nonumber
&
\le 
\frac{9}{\gamma^2}\frac{n}{m^2} x_u^3x_v^3
\langle y^2\rangle
\langle y\rangle^4.
\end{align}
Combining (\ref{2024-10-25+10}),
(\ref{2024-10-25+11})
and
(\ref{2024-10-25+12})
we obtain
\begin{align*}
S_2
&
=
\sum_{\{u,v\}\subset W}\Var |\Lambda_{uv}|
\\
&
\le
n\left(
3
\langle x^2\rangle^2
\langle y^2\rangle
\langle y\rangle^2
+
\frac{9}{\gamma}
\langle x^3\rangle
\langle x^2\rangle
\langle y^2\rangle
\langle y\rangle^3
+
\frac{9}{2\gamma^2}
\langle x^3\rangle^2
\langle y^2\rangle
\langle y\rangle^4
\right)
=
O(n).
\end{align*}

Let us evaluate $S_1$. To this aim we write $S_1$ in the form 
\begin{align*}
S_1=
\sum_{\{u,v\}\subset W}\sum_{t\in W\setminus\{u,v\}}
 \left(
 \Cov (|\Lambda_{uv}|,|\Lambda_{ut}|)
 + 
 \Cov (|\Lambda_{uv}|,|\Lambda_{vt}|)
 \right)
 \end{align*}
and evaluate $ \Cov (|\Lambda_{uv}|,|\Lambda_{ut}|)$
and
$\Cov (|\Lambda_{uv}|,|\Lambda_{vt}|)$.
We now focus on $ \Cov (|\Lambda_{uv}|,|\Lambda_{ut}|)$.
We fix $u,v,t$ and  use  notation (\ref{2024-10-25+30}). In addition, we   denote
\begin{align*}
&
{\bf I}_{(i,j,k)}
=
{\mathbb I}_{iu}
 {\mathbb I}_{ju}
{\mathbb I}_{jt}
{\mathbb I}_{kt},
\qquad
\qquad
\quad
{\bar {\bf I}}_{(i,j,k)}
=
{\bf I}_{(i,j,k)}
-
\E {\bf I}_{(i,j,k)}.
\end{align*}
In view of the second relation of (\ref{2024-10-22}) we have 
 \begin{align*}
 |\Lambda_{uv}|-\E |\Lambda_{uv}|
 =
 \sum_{(i,j,k)\in V_0^3}
 {\bar {\mathbb I}}_{(i,j,k)},
 \qquad
  |\Lambda_{ut}|-\E |\Lambda_{ut}|
  =
  \sum_{(p,q,r)\in V_0^3}
  {\bar {\bf I}}_{(p,q,r)}.
 \end{align*}
 Furthermore, since  for $\{i,j\}\cap \{p,q\}=\emptyset$  random variables 
 $ {\bar {\mathbb I}}_{(i,j,k)}, {\bar {\bf I}}_{(p,q,r)}$ are independent we have
 \begin{align}
 \label{2024-11-02}
 \Cov (|\Lambda_{uv}|,|\Lambda_{ut}|)
&
 =
  \sum_{(i,j,k)\in V_0^3}
  \
  \sum_{(p,q,r)\in V_0^3}
 \E\left(
 {\bar {\mathbb I}}_{(i,j,k)}
  {\bar {\bf I}}_{(p,q,r)}
  \right)
\\
\nonumber
&
=\sum_{(i,j,k)\in V_0^3}
  \sum_{\substack{(p,q,r)\in V_0^3:\\
  \{p,q\}\cap\{i,j\}\not=\emptyset}}
 \E\left(
 {\bar {\mathbb I}}_{(i,j,k)}
  {\bar {\bf I}}_{(p,q,r)}
  \right) 
 \\
 \nonumber
 & 
=:S'_1(u,t;v)+S'_2(u,t;v).
  \end{align}
  Here,  for $h=1,2$, we denote
\begin{align*}
S'_h(u,t;v)
=  \sum_{(i,j,k)\in V_0^3}
  \sum_{\substack{(p,q,r)\in V_0^3:\\
  |\{p,q\}\cap\{i,j\}|=h}}
 \E\left(
 {\bar {\mathbb I}}_{(i,j,k)}
  {\bar {\bf I}}_{(p,q,r)}
  \right).
  \end{align*}
 Now we estimate  $S'_h(u,t;v)$ using inequalities (\ref{covariance}) and $p_{ab}\le p_{ab}^{\star}\le x_by_a/\mu$, see (\ref{2024-09-10}).
  
  Estimation of $S'_2(u,t;v)$. We have
  \begin{align*}
 S'_2(u,t;v)
 &
 =
 \sum_{(i,j,k)\in V_0^3}
 \
 \sum_{r\in V\setminus\{i,j\}} 
 \E\left(
   {\bar {\mathbb I}}_{(i,j,k)}
 ({\bar {\bf I}}_{(i,j,r)}+{\bar {\bf I}}_{(j,i,r)})
 \right)
 \\
 &
 \le
  \sum_{(i,j,k)\in V_0^3}
 \
 \sum_{r\in V\setminus\{i,j\}} 
 \E\left(
  {\mathbb I}_{(i,j,k)}
 ({\bf I}_{(i,j,r)}+{\bf I}_{(j,i,r)})
 \right)
 \\
 &
 =
  \sum_{(i,j,k)\in V_0^3}
 \
 \sum_{r\in V\setminus\{p,q\}} 
p_{iu}p_{ju}p_{jv}p_{kv}
(p_{jt}p_{rt}
+
p_{it}p_{rt}
)
\\
&
\le
\frac{n^4}{\mu^6}x_u^2x_v^2x_t^2
\left(
\langle y^3\rangle
\langle y \rangle^3
+
\langle y^2\rangle^2
\langle y \rangle^2
\right).
 \end{align*}

Estimation of $S'_1$. Denote for short $a_{(i,j,q,r)}
=
{\bf I}_{(i,q,r)}
+   
{\bf I}_{(q,i,r)}
+
{\bf I}_{(j,q,r)}
+   
{\bf I}_{(q,j,r)}
$ 
and ${\bar a}_{(i,j,q,r)}=a_{(i,j,q,r)}-\E \,a_{(i,j,q,r)}$.
We have
  \begin{align*}
 S'_1(u,t;v)
 &
 =
 \sum_{(i,j,k)\in V_0^3}
 \
 \sum_{\substack{(q,r)\in V^2:\, q\not=r,\\
 \{q,r\}\cap \{i,j\}=\emptyset}}
 \E\left(
   {\bar {\mathbb I}}_{(i,j,k)}
 {\bar a}_{(i,j,q,r)}
 \right).
 \end{align*}
 Invoking inequalities
 \begin{align*}
  \E\left(
   {\bar {\mathbb I}}_{(i,j,k)}
 {\bar a}_{(i,j,q,r)}
 \right)
&
 \le
 \E\left(
   {\mathbb I}_{(i,j,k)}
 a_{(i,j,q,r)}
 \right)
 \\
 &
 =
 p_{iu}p_{ju}p_{jv}p_{kv}
 \left(
 p_{qu}p_{qt}p_{rt}
 +
 p_{qu}p_{it}p_{rt}
 +
 p_{qu}p_{jt}p_{rt}
 \right)
 \\
 &
 \le
 \frac{1}{\mu^6}
 x_u^3x_v^2x_t^2 
 \left(
 2
 y_iy_j^2y_ky_q^2y_t
 +
 y_i^2y_j^2y_ky_qy_r
 +
 y_iy_j^3y_ky_qy_r
 \right)
 \end{align*} 
  we obtain
\begin{align*}
   S'_1(u,t;v)
   \le
   \frac{n^5}{\mu^7}
    x_u^3x_v^2x_t^2 
\left(
3
 \langle y^2\rangle^2
 \langle y\rangle^3
+
 \langle y^3\rangle
 \langle y\rangle^4
 \right).
 \end{align*}
  Combining the upper bounds for $S'_1 (u,v;t)$, $S'_2 (u,v;t)$ above  with
  (\ref{2024-11-02}) we obtain
  \begin{align*}
  \sum_{\{u,v\}\subset W}\sum_{t\in W\setminus\{u,v\}}
 \Cov (|\Lambda_{uv}|,|\Lambda_{ut}|)
 &
 \le
 n\langle x^2\rangle^3
 \left(
  \langle y^3\rangle
\langle y \rangle^3
+
\langle y^2\rangle^2
\langle y \rangle^2
\right)
\\
&
+
\frac{n}{\gamma}
\langle x^3\rangle
\langle x^2\rangle^2
\left( 3
 \langle y^2\rangle^2
 \langle y\rangle^3
+
 \langle y^3\rangle
 \langle y\rangle^4
 \right).
 \end{align*}
 The same  bound holds for 
 $\sum_{\{u,v\}\subset W}\sum_{t\in W\setminus\{u,v\}}
 \Cov (|\Lambda_{uv}|,|\Lambda_{vt}|)$. Hence $S_1=O(n)$.
 Proof of Theorem \ref{T3+} is complete.

\end{document}